\documentclass[11pt]{amsart}

\usepackage{amssymb,amsthm,amsmath}
\usepackage[numbers,sort&compress]{natbib}
\usepackage[margin=1in]{geometry}
\usepackage{color}
\usepackage{graphicx}
\usepackage{tikz}
\usepackage[colorlinks,
linkcolor=red,
anchorcolor=blue,
citecolor=blue
]{hyperref}

\def\e{\varepsilon}

\let\newpf\proof \let\proof\relax 
\newenvironment{pf}{\newpf[\proofname]}{\qed\endtrivlist}

\newcommand{\ba}{\overline{A}}

\def\be{\begin{equation}}
\def\ee{\end{equation}}

\def\ba{{\begin{align}}}
\def\ea{{\end{align}}}

\def\bm{\begin{matrix}}
\def\em{\end{matrix}}

\def\0{{\mathbf 0}}

\newtheorem{Theorem}{Theorem}[section]
\newtheorem{Lemma}{Lemma}[section]
\newtheorem{Proposition}{Proposition}[section]
\newtheorem{Corollary}{Corollary}[section]
\newtheorem{Remark}{Remark}[section]

\numberwithin{equation}{section}

\theoremstyle{definition}

\def\tr{{\text{tr}}}

\newcommand{\C}{{\mathbb C}}

\newcommand{\N}{{\mathbb N}}
\newcommand{\Q}{{\mathbb Q}}
\newcommand{\R}{{\mathbb R}}
\newcommand{\T}{{\mathbb T}}

\newcommand{\Z}{{\mathbb Z}}

\def\B0{{\bold{0}}}

\catcode`\@=12

\def\Empty{}
\newcommand\oplabel[1]{
  \def\OpArg{#1} \ifx \OpArg\Empty {} \else
    \label{#1}
  \fi}

\newcommand{\comm}[1]{}
\newcommand{\comment}[1]{}

\begin{document}

\title[]{Multiplicative Jensen's formula and quantitative global theory of one-frequency Schr\"odinger operators}

\author {Lingrui Ge}
\address{Beijing International Center for Mathematical Research, Peking University, Beijing, China.
	} \email{gelingrui@bicmr.pku.edu.cn}

\author {Svetlana Jitomirskaya}
\address{
Department of Mathematics, University of California, Irvine CA, 92697 
} \email{szhitomi@uci.edu}

\author{Jiangong You}
\address{
	Chern Institute of Mathematics and LPMC, Nankai University, Tianjin 300071, China} \email{jyou@nju.edu.cn}

\author{Qi Zhou}
\address{Chern Institute of Mathematics and LPMC, Nankai University, Tianjin 300071, China
} \email{qizhou@nankai.edu.cn}

\begin{abstract}
We introduce the concept of dual Lyapunov exponents, leading to a
multiplicative version of the classical Jensen's formula for
one-frequency analytic Schr\"odinger
cocycles. This formula, in particular, gives a new proof and a
quantitative version of the fundamentals of Avila's global
theory \cite{avila}, fully explaining the behavior of complexified
Lyapunov exponent through the dynamics of the dual cocycle. The key concepts of
(sub/super) critical regimes and acceleration are all explained (in a
quantitative way) through the duality approach. In particular, for
trigonometric polynomial potentials, we establish partial
hyperbolicity of the dual symplectic cocycle and show that the acceleration is
equal to half the dimension of its center, this holding also in the
appropriate sense for the general analytic case. These results lead to a number
of powerful spectral and physics applications.
\end{abstract}

\maketitle
\tableofcontents

\section{Introduction}
\subsection{Multiplicative Jensen's formula}

 Let $f(z)$ be an analytic function given by
 $f(z)=\sum_k\hat{f}(k)z^k$  in $D:=\{x:|z|<r\}$. Suppose that $z_1$, $z_2$, $\cdots$, $z_n$ are the zeros of $f$ in the interior of $D$ repeated according to multiplicity.

 The classical
 {\bf Jensen's formula}, says that for any $0 \leq \e< \ln r$,
\begin{equation}\label{jensen}
   I_\e(f):=\frac{1}{2\pi}\int_0^{2\pi}\ln |f(e^{\e}e^{ix})|dx =I_0(f)-\sum_{\{i:0\leq \ln |z_i|<{\e}\}} \ln |z_i|+\#\{i:0\leq  \ln |z_i|<\e\})\e.
\end{equation}
Using the ergodic theorem, the logarithmic integral on the left hand side
can be interpreted  dynamically, as the limit of time averages along
the trajectory of an ergodic dynamical system.  In particular, given any irrational $\alpha$, one can rewrite \eqref{jensen} as
\begin{align}\label{jensenergodic}
\nonumber & \lim\limits_{n\rightarrow\infty}\frac{1}{2\pi n} \int_0^{2\pi} \ln|f(e^\e e^{i(x+(n-1)\alpha)})\cdots f(e^\e e^{ix})|dx\\
=&I_0(f)- \sum_{\{i:0\leq \ln |z_i|<{\e}\}} \ln |z_i|+\#\{i:0\leq  \ln |z_i|<\e\})\e.
\end{align}
The left hand side of \eqref{jensenergodic} can now be further
interpreted as the complexified Lyapunov exponent of an analytic quasiperiodic  $SL(1,\C)$ cocycle
$(\alpha,f):\T\times\C\to \T\times\C$ that acts via $(\alpha,f)(x,v)=(x+\alpha,f(x)v).$

It is then natural to ask whether there is an analogous formula
for the Lyapunov exponents of matrix-valued cocycles $(\alpha,A)$ where $A$
is an analytic matrix, the situation that is of course much more complicated since the commutativity is lost. The most
intriguing question in this regard is what plays the role of zeros of
analytic function $f$ in the matrix-valued case.


In this paper, we establish such formula for analytic Schr\"odinger
cocycles. In reference to the relation between Birkhoff ergodic
theorem and Kingman's multiplicative ergodic theorem, we call it
multiplicative Jensen's formula.

Schr\"odinger cocycles play a central role in the analysis of one
dimensional discrete ergodic Schr\"odinger operators, a topic with
origins in and a strong ongoing connection to physics and significant
exciting recent advances, particularly in the analytic one-frequency
quasiperiodic case.

Let $\alpha\in \R\backslash\Q$,  $x\in\R,$ and $V$ be a $1$-periodic real analytic
function which can be analytically extended to the strip $\{z||\Im
z|<h\}$. A one-dimensional quasiperiodic Schr\"odinger operator
$H_{V,x,\alpha}:\ell^2(\Z)\to \ell^2(\Z) $ with
one-frequency analytic potential is given by
\begin{equation}\label{so}
(H_{V,x,\alpha}u)_n=u_{n+1}+u_{n-1}+V(x+n\alpha)u_n,
\end{equation}
   The corresponding family of Schr\"odinger cocycles
   $(\alpha,A_E):\T\times\C^2\to \T\times\C^2,\; E\in\R$ is
defined by $(\alpha,A_E)(x,v)=(x+\alpha,A_E(x)v)$ where  $$A_E(x)=\begin{pmatrix}E-V(x)&-1\\ 1&0\end{pmatrix}.$$ It governs the
behavior of solutions to $$\label{eigen} H_{V,x,\alpha}u=Eu.$$
The complexified Lyapunov exponent is given by
\begin{equation}\label{multiergodicsch}
L_\e(E)=\lim\limits_{n\rightarrow\infty}\frac{1}{n} \int \ln \|A(x+i\e+(n-1)\alpha)\cdots A(x+i\e)\|dx.
\end{equation}
The limit existing, as usual, by the Kingman's subadditive ergodic
theorem. Complexified Lyapunov exponents were first studied by
M. Herman \cite{her}, were crucial in the proofs of positivity
of Lyapunov exponents at large couplings \cite{ss,bg,bourg} and
played a central role in Avila's global theory \cite{avila}.

We establish an analogue of \eqref{jensenergodic} for $L_\e(E),$
where it turns out that the role of zeros of $f$ in \eqref{jensen} is
played by the (appropriate limits of) the
Lyapunov exponents of the dual cocycles, an object that we prove to exist and call {\bf
  dual Lyapunov exponents}.

The Aubry dual of the one-frequency Schr\"odinger operator \eqref{so} is
\begin{equation}\label{long}
(\widehat{H}_{V,\theta,\alpha}u)_n=\sum\limits_{k=-\infty}^{\infty} V_k u_{n+k}+2\cos2\pi(\theta+n\alpha)u_n, \ \ n\in\Z.
\end{equation}
where $V_k$ is the $k$-th Fourier coefficient of $V,$ see Sec \ref{strip} for details.
For general analytic $V,$ operator \eqref{long} is infinite-range, so
its eigenequation $\widehat{H}_{V,\theta,\alpha}u=Eu$ does not define any
cocycle. However, if $V(x)$ is a
trigonometric polynomial of degree $d,$ the eigenequation
$\widehat{H}_{V,\theta,\alpha}u=Eu$ leads to a symplectic $2d$-dimensional cocycle
that we denote  by
$(\alpha,\widehat{A}_E)$. We denote its Lyapunov exponents
by $\pm\hat{L}^d_{1}(E), \cdots, \pm\hat{L}^d_{d}(E)$ according to
multiplicity \footnote{See Section \ref{com} for the definitions and
  discussion.}. We may assume $0\leq \hat{L}^d_1(E)\leq \cdots\leq
\hat{L}^d_d(E)$. We have

\begin{Theorem}\label{1}
Assume $V(x)$ is a trigonometric polynomial of degree $d.$ For
$\alpha\in\R\backslash\Q$ and $(E,\e)\in\R^2$, we have
\begin{align*}\label{gne1}
L_{\e}(E)= L_0(E) -\sum_{\{i:\hat{L}^d_i(E)< 2\pi |\e|\}}  \hat{L}^d_i(E)+2\pi(\#\{i:\hat{L}^d_i(E)<2\pi|\e|\})|\e|.
\end{align*}
\end{Theorem}

In fact, the multiplicative Jensen's formula \eqref{1} is not merely
an analogue of the classical Jensen's formula but a proper
generalization, because zeros of an analytic function $f$ can also be
interpreted as the Lyapunov exponents of the dual cocycle.
 Indeed, consider the diagonal operator acting on $\ell^2(\Z)$
\begin{equation}\label{multiplication1}
(M_{x} u)_n=V(x+n\alpha)u_n, \ \ n\in\Z,
\end{equation}
where $V$ is a $1$-periodic real trigonometric polynomial of degree $d.$
Its Aubry dual is given by the T\"oplitz operator
\begin{equation}\label{multiplicationdual1}
(\widehat{M} u)(n)=\sum\limits_{k=-d}^dV_ku_{n+k}, \ \ n\in\Z.
\end{equation}

 It turns out that if $\{z_1(E),\cdots, z_d(E)\}$ are zeros of $V(z)=E$ with
  $1\leq |z_i(E)|,$ \footnote{The zeros come in pairs because $V$ is
   real.} then $\pm\ln|z_i|$ are precisely the Lyapunov exponents of
   the  cocycle $(\alpha,\widehat{M})$    \footnote{Since $(\alpha,\widehat{M})$ is a constant cocycle, its
     Lyapunov can be easily calculated.} of the eigenequation $\widehat{M}u=Eu,$
   while $I_{\e}(E):=\int_{0}^{1}\ln|E-V(x+i|\e|)|dx$ is the
complexified Lyapunov exponent  of the  $SL(1,\C)$ cocycle
$(\alpha,V):\T\times\C\to \T\times\C$ acting via $(\alpha,V)(x,v)=(x+\alpha,V(x)v)$.






If $V$ has infinitely many harmonics, we will use trigonometric
polynomial approximation. Let $V^d(x)=
\sum_{k=-d}^d \hat{V}_ke^{2\pi i kx}$ and let $\hat{L}^d_i(E)$ be the
Lyapunov exponents of the corresponding dual ${\rm Sp_{2d}(\C)}$ cocycle. We have

\begin{Theorem}\label{1general}{\bf [The multiplicative Jensen’s formula]}
For $\alpha\in\R\backslash\Q$ and $V\in C^\omega_h(\T,\R)$, there exist  non-negative $\{\hat{L}_i(E)\}_{i=1}^m$ such that for any $E\in\R$
$$
\hat{L}_i(E)=\lim\limits_{d\rightarrow \infty}\hat{L}^d_i(E), \ \ 1\leq i\leq m.
$$
Moreover, \begin{align*}\label{gne1}
L_{\e}(E)= L_0(E) -\sum_{\{i:\hat{L}_i(E)< 2\pi|\e|\}}  \hat{L}_i(E)+2\pi(\#\{i:\hat{L}_i(E)<2\pi|\e|\})|\e|
\end{align*} for $|\e|<h$.
\end{Theorem}
\begin{Remark}Note that the cocycle itself changes dramatically when
  $d$ changes, with no limit in any of its components,
  however the Lyapunov exponents do converge to their limits, that we
  call {\bf dual Lyapunov exponents} of \eqref{so}.
  \end{Remark}
\begin{Remark}
One of the fundamental results in \cite{avila} is  that $L_\e(E)$ is a
piecewise affine function in $\e$ for each $E$, and the slope of each
piece is an integer.  Theorem \ref{1general} quantifies this result,
identifying the turning points with {\it distinct} $\hat{L}_i$'s, and
the increments in the integer slopes with multiplicities of distinct $\hat{L}_i$'s.
\end{Remark}
Indeed, for a fixed $E\in\R$, assume that
$$
0\leq \hat{L}_{k_1}<\hat{L}_{k_2}<\cdots<\hat{L}_{k_{\ell}}
$$
and the multiplicity of each $\hat{L}_{k_i}$ is $\{k_{i}-k_{i-1}\}_{i=1}^\ell$ with $k_0=0$ and $k_\ell=m$. One may rewrite $L_{\e}(E)$ in Theorem \ref{1general} as the following piecewise affine function,
\begin{align}\small
L_{\e}(E)=\left\{
\begin{aligned} & L_0(E) &|\e|\in
\left[0,\frac{\hat{L}_{k_1}}{2\pi}\right],\\
&L_{\frac{\hat{L}_{k_i}}{2\pi}}(E)+2\pi k_{i}\left(|\e|-\frac{\hat{L}_{k_i}(E)}{2\pi}\right)
&|\e|\in \left(\frac{\hat{L}_{k_{i}}}{2\pi},\frac{\hat{L}_{k_{i+1}}}{2\pi}\right],\\
&L_{\frac{\hat{L}_{k_\ell}}{2\pi}}(E)+2\pi k_\ell\left(|\e|-\frac{\hat{L}_{k_\ell}(E)}{2\pi}\right) &|\e|\in \left(\frac{\hat{L}_{k_\ell}}{2\pi},h\right).
\end{aligned}\right.
\end{align}
where $1\leq i\leq \ell-1$. See pictures I-III for three different cases.

\begin{figure}[htbp]
\centering
\begin{minipage}[t]{0.4\textwidth}
\centering
\includegraphics[width=1\linewidth]{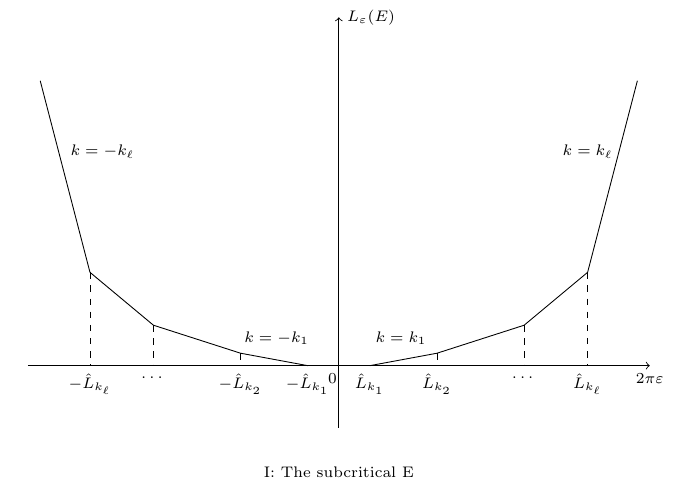}
\end{minipage}
\begin{minipage}[t]{0.4\textwidth}
\centering
\includegraphics[width=1\linewidth]{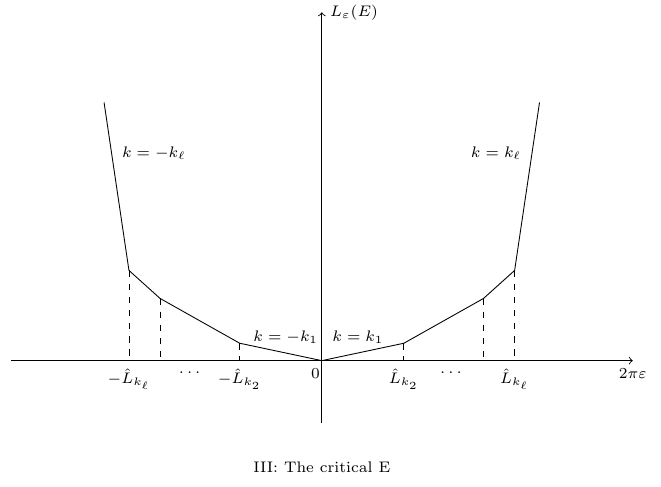}
\end{minipage}
\begin{minipage}[t]{0.4\textwidth}
\centering
\includegraphics[width=1\linewidth]{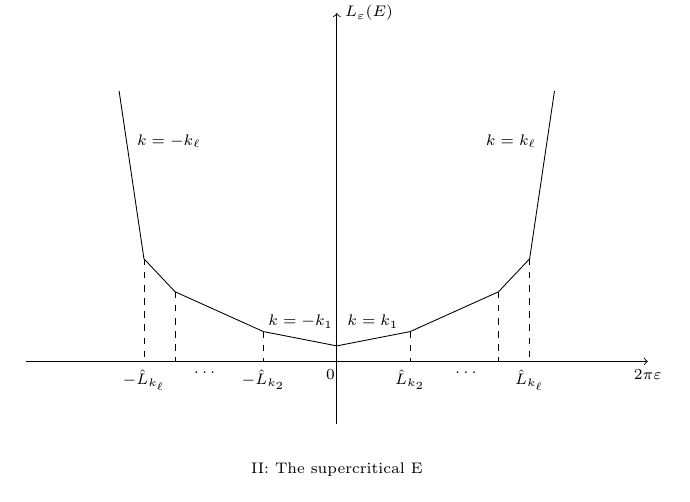}
\end{minipage}
\end{figure}


\subsection{Quantitative version of Avila's global theory}
The multiplicative Jensen's formula not only sheds the light on the global theory of one-frequency Schr\"odinger
cocycles, but allows crucial advances in the study of the spectral
theory of one-frequency Schr\"odinger operators \eqref{so}. 

In the past 40 years after the groundbreaking paper \cite{ds} the theory of quasiperiodic Schr\"odinger
operators has been developed extensively, see
\cite{b1,damanik,jcong,jm,you} for surveys of more recent results. For
the one-frequency case, starting with \cite{SJ94} and
  then \cite{J,bg} the main thread
has been to
establish resuts non-perturbatively, i.e. either in the regime of almost reducibility \cite{puig1,puig2,aj1,afk,hy,avila1,avila2, you} or in the regime
of positive Lyapunov exponent \cite{J,bg,b1,bj02,gs1,gs2,gs3,aj,JLiu,JLiu1}. In 2015,  Avila \cite{avila} gave a qualitative spectral
picture, covering both regimes, 
based on the analysis of the asymptotic behavior of $L_\e(E)$. The
central concept in Avila's global theory \cite{avila} is the {\it
  acceleration}
$$
\omega(E)=\lim\limits_{\e\rightarrow 0^+}\frac{L_\e(E)-L_0(E)}{2\pi\e}.
$$
The global theory divided  the spectra of  one-frequency Schr\"odinger operator  into three regimes based on the Lyapunov exponent and acceleration:
\begin{enumerate}
\item The subcritical regime: $L(E)=0$ and $\omega(E)=0$.
\item The critical regime: $L(E)=0$ and $\omega(E)>0$.
\item The supercritical regime: $L(E)>0$ and $\omega(E)>0$.
\end{enumerate}
Moreover, the subcritical regime is equivalent to the almost reducible
regime \cite{avila1,avila2}. The critical regime is rare in the sense
that it is a set of zero Lebesgue measure  \cite{avila,ak,jkr}. We will use
the (sub/super)critical terminology both when referring to the energies
$E$ and to the corresponding cocycles $(\alpha,A_E).$

Note that the global theory terminology was motivated by the study of the almost Mathieu operator (AMO), the central model in one-frequency quasiperiodic Schr\"odinger operators,
\begin{equation}\label{amo}
(H_{\lambda,x,\alpha}u)_n=u_{n+1}+u_{n-1}+2\lambda\cos2\pi(x+n\alpha)u_n,\ \ n\in\Z,
\end{equation}
where explicit computation \cite{bj,avila} shows that for all $E$ in the spectrum, we have
\begin{enumerate}
\item $|\lambda|<1$: $L(E)=0$ and $\omega(E)=0$.
\item $|\lambda|=1$: $L(E)=0$ and $\omega(E)=1$.
\item $|\lambda|>1$: $L(E)=\ln|\lambda|$ and $\omega(E)=1$.
\end{enumerate}

Roughly speaking, Avila's global theory is based on the picture of
$L_\e(E)$ for $\e$ small enough. Our multiplicative Jensen's formula
actually not only gives the full picture of $L_\e(E)$ for any
$|\e|<h$, but also gives quantitative characterizations of several
quantities in \cite{avila}, such as the {\it acceleration} and the
{\it subcritical radius} defined below.

In particular, we can recharacterize Avila's (sub/super)critical regimes  in terms of  the Lyapunov
exponents $L(E)$ \footnote{We sometimes identify $L_0(E)$ and $L(E)$.}
and the smallest non-negative ``dual Lyapunov exponent", without using
the concept of acceleration:
\begin{Theorem}\label{2general}
Assume $\alpha\in \R\backslash\Q$ and $V\in C_h^\omega(\T,\R)$, then  $E\in\R$ is
\begin{enumerate}
\item{\it {Outside the spectrum} \footnote{I.e. uniformly
      hyperbolic. See Section \ref{secds} for the definition of uniform hyperbolicity.} if $L(E)>0$ and $\hat{L}_1(E)>0$,}
\item{\it{Supercritical} if $L(E)>0$ and $\hat{L}_1(E)=0$,}
\item{\it{Critical} if $L(E)=0$ and $\hat{L}_1(E)=0$,}
\item{\it{Subcritical} if $L(E)=0$ and $\hat{L}_1(E)>0$.}
\end{enumerate}
\end{Theorem}
\begin{Remark}
(4) implies that the Schr\"odinger cocycle $(\alpha,A_E)$ is subcritical if and only if its ``dual  Lyapunov exponents" are all positive, which serves as the  basis for the first author's new proof of the almost reducibility conjecture \cite{ge}.
\end{Remark}

We also give a new quantitative characterization of Avila's acceleration:
\begin{Corollary}\label{3general}
For $\alpha\in \R\backslash\Q,$  $V\in C_h^\omega(\T,\R)$, and any
$E\in\R$ we have
$$
\omega(E)=\begin{cases}
0 &\hat{L}_1(E)>0\\
\#\left\{j| \hat{L}_j(E)=0\right\} &\hat{L}_1(E)=0
\end{cases}.
$$
\end{Corollary}

\begin{Remark}
The  acceleration plays a crucial role in the study of supercritical
Schr\"odinger operators. Corollary \ref{3general} shows that it is
equal to the number of dual Lyapunov exponents that are equal to zero,
or, for trigonometric polynomial $V$, to the dimension of the center of the corresponding
cocycle. Generally speaking, although the definition of acceleration
is neat, it's not easy to see why the acceleration is an integer
outside the uniformly hyperbolic regime where it is simply equal to
the winding number of a certain function. It is also difficult to compute  the acceleration for specific cocycles. Corollary \ref{3general} provides another point of view which is more convenient, at least in the perturbative case (see Section \ref{ar-ac} for further discussion). 
\end{Remark}

In the study of subscritical Schr\"odinger operator, an important quantity is the so-called {\it subcritical radius} defined by
$$
h(E)=\sup\{|\e|:L_\e(E)=0\}.
$$
It turns out it is also linked to dual
Lyapunov exponents.
\begin{Corollary}\label{4general}
For all $\alpha\in \R\backslash\Q$, $V\in C_h^\omega(\T,\R),$ and
$E\in\R$, we have $h(E)=\frac{\hat{L}_1(E)}{2\pi}$.
\end{Corollary}
\begin{Remark}
For subcritical almost Mathieu operators, it's explicitly computed in \cite{avila} that
$$
h(E)=\frac{\hat{L}_1(E)}{2\pi}=-\frac{\ln|\lambda|}{2\pi}
$$ for all $E$ in the spectrum, which plays an important role in several optimal estimates \cite{gyz,gyz1}. Corollary \ref{4general} is a generalization of this fact to general one-frequency Schr\"odinger operators.
\end{Remark}
\subsection{Aubry duality}
Our work can be viewed in a nutshell as the duality approach to Avila's global theory.
Aubry duality: a Fourier-type transform that links the direct integral
in $x$ of operators \eqref{so} to the direct integral in $\theta$ of
operators \eqref{long} has had a long history since its original
discovery by Aubry-Andre \cite{aa} and has been explored and applied
at many levels. Representing a certain gauge invariance of the
underlying two-dimensional discrete operator in a perpendicular
magnetic field \cite{mandelshtamzhitomirskaya}, it has been understood
at the level of integrated density of states, Lyapunov exponents,
individual eigenfunctions and dynamics of individual cocycles, and
explored in various qualitative and quantitative ways.

The almost Mathieu family stands out among other quasiperiodic
operators \eqref{so} precisely because it is self-dual with respect to the
Aubry duality, with
$\hat{H}_{\lambda,x,\alpha}=\lambda H_{\frac{1}{\lambda},x,\alpha},$ e.g
\cite{as}. In
particular,  the subcritical regime ($|\lambda|<1$) and the
supercritical regime ($|\lambda|>1$) are dual to each other, and this
has been fruitfully explored in both directions.   Aubry duality 
enables one to use the supercritical techniques (localization method)
to deal with the subcritical problems \cite{J,puig1,puig2,aj,aj1,gk},
as well as the subcritical methods (almost reducibility) to study the
supercritical problems \cite{ayz,ayz1,gyz,gyzh1,gy, gyz1,you}. Even though
the self-duality is destroyed when going beyond the almost
Mathieu operator, many of the sub(super)critical results for the
almost Mathieu operator can be generalized to \eqref{so}
or \eqref{long}. Based on the localization method for operator
\eqref{long}, one can get (almost) reducibility results for operators
\eqref{so}, see \cite{bj02,puig1,aj1}. Almost reducibility for
operator \eqref{so} in turn implies localization results for
operator \eqref{long}, see \cite{ayz,jk,gyz,gyz1,gy}.  Aubry duality
therefore serves as a powerful bridge between  \eqref{so}  and
\eqref{long}.

All these methods and connections so far stayed firmly on the real
territory, where both the operator and its dual are self-adjoint, so
one can enjoy all the benefits of the spectral theory.  Here we, for the first time, find the way to complexify the
Aubry duality, or, alternatively, extend it to the {\it
  non-self-adjoint} setting, leading both to a new manifestation of it
and a new empirical understanding, as well as a much
deeper understanding of the existing manifestations.

 Historically, Aubry duality was first formulated at the level
of the integrated density of states, and thus, using the Thouless
formula, the Lyapunov exponents. Namely, it was shown in \cite{aa} (with the
argument made rigorous in \cite{as}) that for the almost
Mathieu operator $H_{\lambda,x.\alpha}$ given by \eqref{amo}, the
following relation holds
\begin{equation}\label{aaas} L(E) = \hat{L}(E)+\log|\lambda|\end{equation}

A similar argument based on the Thouless formula for the strip
\cite{ks} leads to the beautiful Haro-Puig formula \cite{haro} for
operators \eqref{so} with
trigonometric polynomial $V(x)$
\begin{equation}\label{hpu}
L(E)=\sum_{\{i:\hat{L}^d_i(E)>0\}}  \hat{L}^d_i(E) +\ln |V_{d}|,
\end{equation}
which specializes to \eqref{aaas} when $V(x)=2\cos 2\pi x,$ since in
this case $V_{d} =V_{-d}=\lambda.$ The
multiplicative Jensen's formula that we discover can be manipulated
into the Haro-Puig formula \eqref{hpu} for complexified Lyapunov
exponents, but the latter cannot be seen in the framework of
the existing proof, in absence of self-adjointness and the related
spectral theory based invariance of the integrated density of states,
and in itself presents no compelling reason for it to hold.

We discover, however, that Aubry duality can be understood in a way
that does not require any self-adjointness, leading to
a new dynamical  perspective  on it and playing an important role in
enabling  various spectral applications.

We show that the fundamental way to see Aubry duality is through the
invariance of the averaged Green's function 
$$
 \int_{\T}\langle \delta_0,(H_{V(\cdot+i\e),x,\alpha}-E)^{-1}\delta_0\rangle dx = \int_{\T}\langle \delta_0, (\hat{H}_{V(\cdot+i\e),\theta,\alpha}-EI)^{-1}\delta_0\rangle d\theta,
 $$
 something that can then be approached dynamically and combined with a non-self-adjoint version
 of the Johnson-Moser's theorem \cite{JM} that links averaged Green's
 function to the {\it derivative} of the Lyapunov exponent, a strategy that we discuss
 more in Section \ref{ideas}.

 The classical empirical understanding of Aubry duality is that
 Fourier transform takes nice normalizable eigenfunctions into Bloch
 waves and vice versa. Alternatively, (almost) localized
 eigenfunctions correspond to (almost) reducibility for the dual
 cocycle, and vice versa, something that by now has almost became a
 folklore. Here we present a similarly complelling heuristics - a new
 perspective - that was behind our discovery of the multiplicative Jensen's formula.
 
Assume that  $(\alpha, \widehat{A}_{E})$ is analytically conjugated to the form:
\begin{equation}\label{co-au} Z(\theta+\alpha)^{-1} \widehat{A}_{E}(\theta)Z(\theta) = \left(
  \begin{array}{ccccccc}
e^{\gamma} & 0& \quad  \\
0 & e^{-\gamma} & \quad  \\
\quad & \quad &  D(\theta) \end{array}
  \right) .\end{equation}
By the Aubry duality, it implies that 
\begin{equation*}
(H_{V,x+i\gamma,\alpha}u)_n=u_{n+1}+u_{n-1}+V(x+i\gamma+n\alpha)u_n,
\end{equation*}
has a localized eigenfunction. Therefore the Schr\"odinger cocycle $(\alpha, 
 A_E(x+i\gamma))
  $ is nonuniformly hyperbolic, so $L_{\e}(E)$ cannot be affine at
  $\e=\gamma,$ therefore $\gamma$ must be the turning point of
  $L_{\e}(E)$. But of course \eqref{co-au} just means $\gamma$ is the Lyapunov exponent of the dual cocycle  $(\alpha, \widehat{A}_{E})$.
While not fully rigorous, we see this argument as inspirational to our
approach, and in fact it plays an important role both in the final
proof and a physics application \cite{Liu2020L,LZC}.

\subsection{Bochi-Viana Theorem for dual cocycles and
partial hyperbolicity}

Both our proof of Theorem \ref{1} and
an important starting point for the most interesting corollaries
is based on the study of the dynamics of dual cocycles,
which turn out to have a remarkable universal property.

It is a general program, first outlined by  Ma\~n\'e \cite{m} and
developed by Bochi-Viana \cite{bv} that,  when applied to linear
cocycles, states that for $C^0$ generic ${\rm GL_d(\C)}$ cocycles over any
measure preserving transformation the Oseledets
splitting (see section \ref{secds} for the definitions in this setting)
is either trivial or dominated. While this result definitely hinges
on low regularity considerations (and counterexamples in higher
regularity do exist), it was shown in \cite{ajs} that Bochi-Viana
theorem also holds - and in a much stronger form - for analytic
one-frequency cocycles: the Oseledec splitting is either trivial or
dominated on an open and dense set of such cocycles.

Here we show that something stronger yet holds for the dual
cocycles. Let $\C_+$  denote $\{E\in\C|\Im E>0\}.$ For $V(x)$  a trigonometric polynomial of degree
  $d,$ let $$
0\leq \hat{L}_{k_1}<\hat{L}_{k_2}<\cdots<\hat{L}_{k_{\ell}}
$$ be 
the listing of all nonnegative dual Lyapunov esponents, where the multiplicity of each $\hat{L}_{k_i}$ is $\{k_{i}-k_{i-1}\}_{i=1}^\ell$ with $k_0=0$ and $k_\ell=d$.
\begin{Theorem} \label{dom} Let $V(x)$ be a trigonometric polynomial of degree
  $d.$ Then the dual 
  cocycle
  $(\alpha,\widehat{A}_{E})$ is {\it always}
  \begin{enumerate} \item
$(d-k_i)$-dominated for all
      $ 0\leq i \leq \ell,$
       for $E\in\C_+$;
     \item either trivial or
$(d-k_i)$-dominated for all
      $1\leq i \leq \ell,$
       for $E\in\R.$
    \end{enumerate}
    \end{Theorem}
    \begin{Remark} For $E\in\C_+$ the cocycle is obviously uniformly
      hyperbolic, so $d$-dominated, but the domination in all other $k$
      is a nontrivial statement.
    \end{Remark}
In particular, we have
\begin{Corollary}\label{ph} The acceleration $\omega(E)>0$ if and only
  if the dual ${\rm Sp_{2d}(\C)}$ cocycle
  $(\alpha,\widehat{A}_{E})$ is partially hyperbolic 
    with zero  center Lyapunov exponents.
  \end{Corollary}

\subsection{A spectral application}
In this subsection, we give a sample direct spectral application of
our quantitative global theory: a new neat
characterization of the spectrum of $H_{V,x,\alpha}$
and a criterion for uniformity of corresponding Schr\"odinger cocycles.

It is well-known
that the spectrum of $H_{V,x,\alpha},$ denoted $\Sigma_{V,\alpha},$ is
an $x$-independent set \cite{as}.  The classical Johnson's theorem \cite{johnson} characterizes the spectrum as
$$
\Sigma_{V,\alpha}=\{E\in\R: L(E)=0\ \ \mbox{or}\ \  
(\alpha, A_E)\ \ \mbox{is
 non uniformly hyperbolic}\}.
$$

Non-uniform hyperbolicity is generally difficult to capture. It turns
out however that it is determined precisely by the lowest dual
Lyapunov exponent. We have

\begin{Corollary}\label{spec cha}
For any $\alpha\in \R\backslash\Q$ and $V\in C_h^\omega(\T,\R)$, then
$$
\Sigma_{V,\alpha}=\{E\in\R: L(E)\cdot \hat{L}_1(E)=0\}.
$$
\end{Corollary}

An equivalent formulation of Corollary \ref{spec cha} is the following
criterion for uniformity of 
Schr\"odinger cocycles.
 We recall that an $SL(2,\C)$ cocycle  $(\alpha,A)$ is  {\it
  uniform} if the convergence 
$$
\lim\limits_{n\rightarrow \infty }\frac{1}{n}\ln\|A(x+(n-1)\alpha)\cdots A(x)\|=L(\alpha,A)
$$
holds for {\it all}
$x\in\T$ and is uniform (see, e.g. \cite{DL} for a discussion).
Since Schr\"odinger cocycles $(\alpha,A_E)$ are uniform for $E$
outside the spectrum or in the set where $L(E)=0$
(e.g. \cite[Corollary A.3]{DL}), an immediate consequence of Corollary \ref{spec cha} is

\begin{Corollary}\label{uniform}
 A Schr\"odinger cocycle $(\alpha,A_E)$ with  $\hat{L}_1(E)>0$ is
 always uniform.
\end{Corollary}

\begin{Remark} If $V$ is a trigonometric polynomial, this can be
  neatly reformulated as `` A Schr\"odinger cocycle
    with hyperbolic dual cocycle  is always uniform". 
  \end{Remark}
  
Most excitingly however, our analysis enables to extend some of
the most famous almost Mathieu results to large classes of
quasiperiodic operators.

In particular, in the companion paper \cite{martini} we develop machinery to
prove the {\it Ten Martini problem} (that is Cantor spectrum without any parameter exclusion) for a
large explicitly defined open set of both sub and
super critical quasiperiodic operators, so
called operators of type 1. The Ten Martini problem has so far
only been established for the almost Mathieu operator through a
combination of Liouville and Diophantine approaches that were both
almost Mathieu specific and only quite miraculously
met in the middle. It has not even been universally expected that it
holds for all parameters for anything other than the almost Mathieu
operator.

In \cite{gj} we prove {\it sharp arithmetic spectral transition},
as in \cite{ayz,JLiu,JLiu1} for all operators of type 1, without further
assumptions.

Finally, these results enable a new and simple proof of Avila's almost reducibility
conjecture for Schr\"odinger cocycles \cite{ge}. With
  subciticality guaranteeing $\hat{L}_1>0,$ the
  proof proceeds through establishing non-perturbative almost
  localization for the dual operator and is
  optimal for the case of trigonometric polynomial $V$ (i.e. does not require  shrinking
of the band).

This paper is a result of a long-term effort. Our results, have been obtained
and proofs fully written several years ago, with this release delayed
by our quest to obtain the strongest applications. The latter is an
ongoing and expanding project. The results, as well as
some of the applications, have already
been presented and discussed, with many details, at multiple venues,
including the Anosov-85 meeting, November 2021, BIRS workshop on almost-periodic spectral
    problems, April 2022,
ICM 2022, and QMath 15, 2022, as well as announced in
  \cite{jcong} \footnote{A preprint was made
    available to the community in January 2022.}. Several physics-related applications of
our results
\cite{Liu2020L,LZC}, with \cite{LZC} made rigorous in \cite{WWYZ},
have already 
appeared and the applications to the ten
martini problem \cite{martini} and the almost reducibility conjecture
\cite{ge} are being released simultaneously. \\

\smallskip

The rest of this paper is organized as follows. Section \ref{appl}
contains further spectral and physics applications. In Section
\ref{ideas}, we introduce the main ideas of the proof. Section
\ref{pre} contains the preliminaries. In Section \ref{hermit}, we
study the Green's function of general finite-range
Schr\"odinger operators, while in Section \ref{non-hermit} we study
the Green's function for non-self-adjoint quasiperiodic  operators.
In Section \ref{main-proof}, we prove the main results,  postponing
proofs of the remaining results to Section \ref{remain}.  Finally, we prove Johnson-Moser's theorem for Schr\"odinger operators on the
    strip (Proposition \ref{mg}) in Section \ref{proofmg}, and prove the  representation of the Green’s function for general strip
    operators (Lemma \ref{basic})  in   Section \ref{proofbasic}.

\section{Other applications}\label{appl}

\subsection{Arithmetic Anderson localization}

Our results allow us to make spectral conclusions both for
$H_{V,x,\alpha}$ and $\widehat{H}_{V,\theta,\alpha}$.  Here we present a
sample result on Anderson localization for
$\widehat{H}_{V,\theta,\alpha}$,  which
 was extensively studied
\cite{aj1,gyz,bj,cd,jk,b1,gy} since the 1980s. All  the existing
results are ``local" in the sense that one needs to assume there is a
large coupling constant $\lambda$ before the $\cos$
potential. moreover most of the results cannot go beyond the Diophantine frequencies. We give a global result, starting from the positivity of the Lyapunov exponents. Let
$$
\beta=\beta(\alpha)=\limsup\limits_{k\rightarrow \infty}-\frac{\ln\|k\alpha\|_{\R/\Z}}{|k|}.
$$
For a given irrational number $\alpha$, we say $\theta\in (0,1)$ is $\alpha$-Diophantine if there exist $\kappa>0$ and $\tau>1$ such that
$$
\|2\theta+k\alpha\|_{\R/\Z}>\frac{\kappa}{(|k|+1)^\tau},
$$
for any $k\in\Z$, where $\|x\|_{\R/\Z}=\text{dist}(x,\Z).$ Clearly, for any fixed irrational number $\alpha$, the set of phases which are  $\alpha$-Diophantine  is of full Lebesgue measure.

\begin{Corollary}\label{localization}
If $\hat{L}_1(E)>\beta>0$ for all $E\in\R$, then $\widehat{H}_{V,\alpha,\theta}$ has Anderson localization for $\alpha$-Diophantine $\theta$.
\end{Corollary}
\begin{Remark}
For the almost Mathieu operator, Corollary \ref{localization} is what
is now sometimes called the Andre-Aubry-Jitomirskaya conjecture
\cite{aa,j94} which was proved in \cite{JLiu}, see also \cite{gyzh1} for a new proof.
\end{Remark}
\begin{Remark}The limitation $\beta>0$ comes from our reliance in the
  proof on a theorem of \cite{gyzh1}, who in turn rely on Avila's
   proof of the almost reducibility conjecture for Liouville
   $\alpha$ \cite{avila1}. This limitation has been removed in the
   follow-up paper by the first author \cite{ge} through a direct
   localization-side proof for $\beta=0.$
   \end{Remark}
   \begin{Remark} We present the result for $\alpha$-Diophantine
     $\theta$ rather than a slightly weaker optimal \cite{JLiu}
     condition $\delta(\alpha,\theta)=0$ where
     $$\delta(\alpha,\theta)=\limsup_{k\to\infty} -\frac{\ln ||2\theta+
       k\alpha||_{\R/\Z}}{|k|},$$ because the authors of \cite{gyzh1}
     choose a similar limitation. The theorem in fact holds under the
     $\delta(\alpha,\theta)=0$ condition with a little more technical effort.
    \end{Remark} 

\subsection{An application to the Soukoulis-Economou's
model}

We can also make immediate conclusions for the Soukoulis-Economou's
model (SEM)
\begin{equation}\label{gamo}
(H_{\alpha,x}u)(n)=u(n+1)+u(n-1)+2\lambda_1\cos2\pi(x+n\alpha)u(n)+2\lambda_2\cos
4\pi(x+n\alpha)u(n).
\end{equation}
It is also known in physics literature as generalized Harper's model
(e.g. \cite{hk,se}), which is of special interest because of its
connection to the three dimensional quantum Hall effect
\cite{hk,se}. The Lyapunov exponents for this model have been studied
in \cite{jliu3, shoumarx}.

The Aubry dual of (\ref{gamo}) is
\begin{equation}\label{dual_gamo}
(\widehat{H}_{\theta,\alpha}u)(n)=\lambda_2u(n-2)+\lambda_1u(n-1)+\lambda_1u(n+1)+\lambda_2u(n+2)+2\cos2\pi(\theta+n\alpha)u_n, \ \ n\in\Z.
\end{equation}
The operator \eqref{dual_gamo} is a 4-th order difference operator, and we denote the non-negative Lyapunov exponent of the associated cocycle by $\hat{L}_2(E)\geq \hat{L}_1(E)\geq 0$.
\begin{Corollary}\label{ghm1}
For SEM operator with $\alpha\in\R\backslash\Q,$ for any $E\in \R$, $\omega(E)=2$ if and only if $L(E)=\ln|\lambda_2|$ and $|\lambda_2|\geq 1$.
\end{Corollary}

\begin{Corollary}\label{ghm2}
For $\alpha\in\R\backslash\Q$ and $|\lambda_2|< 1,$ the energies in
the spectrum of SEM are in one of the following  three regimes:
\begin{enumerate}
\item Subcrtical regime: $L(E)=0$ and $\omega(E)=0$.
\item Critical regime: $L(E)=0$ and $\omega(E)=1$.
\item Supercritical regime: $L(E)>0$ and $\omega(E)=1$.
\end{enumerate}
\end{Corollary}
\begin{Remark}
In this case, the crucial point is that the acceleration is always no more
than $1,$ which is  also a key feature of the almost Mathieu operator. In
particular, this means that supercritical SEM with
$|\lambda_2|<1$ is of type 1 in the sense of
\cite{martini}, and it makes it possible to generalize many almost
Mathieu results to this case. We note that the
supercritical regime is known to hold under explicit
conditions on $\lambda_1,\lambda_2$ with $|\lambda_2|<
1$ \cite{jliu3} requiring, in particular, $\lambda_1 >100\lambda_2$. Our
analysis of type 1 operators applies to the entire regime  $|\lambda_2|<
1$.
\end{Remark}

\subsection{A further characterization of the acceleration}\label{ar-ac}
Let $V$ be a trigonometric polynomial of degree $d$ such that
$\widehat{A}_{E}$ is almost reducible to some constant matrix
$\tilde{A}$ in the sense that there exists $B_n\in
C^\omega_{r_n}(\T,{\rm GL_{2d}(\C))}$ for some $r_n>0$ such that
$$\|B_n(\theta+\alpha)^{-1}\widehat{A}_{E}(\theta)B_n(\theta)-\tilde{A}\|_{r_n}\rightarrow
0.$$ Note that this assumption is always satisfied for a positive
measure set of $\alpha$ if $V=\lambda f$  and $\lambda$ is
sufficiently small. In this case, the dual Lyapunov exponents can be
computed explicitly, and the multiplicative Jensen's formula takes a
particularly neat form

\begin{Corollary}\label{6}
Suppose that $E\in\R$, $\alpha\in\R\backslash\Q$ and $(\alpha,\widehat{A}_{E})$ is almost reducible to
some constant matrix $\tilde{A}$.  Let $\lambda_1, \cdots, \lambda_d$ be the eigenvalues of  $\tilde{A}$, counting the multiplicity, with $1\leq |\lambda_1|\leq \cdots\leq |\lambda_d|$. Then
\begin{align*}
L_{\e}(E)= L_0(E) -\sum_{\{i:\ln |\lambda_i|< 2\pi |\e|\}}  \ln|\lambda_i|+2\pi (\#\{i:\ln|\lambda_i|<2\pi|\e|\})|\e|.
\end{align*}
\end{Corollary}

\begin{Corollary}\label{7}
Under the assumptions of Corollary \ref{6}, we have
$$
\omega(E)=\begin{cases}
0 &|\lambda_1|>1\\
\#\left\{j| |\lambda_j|=1\right\} &|\lambda_1|=1
\end{cases}.
$$
\end{Corollary}
\begin{Remark}
The acceleration is nothing but the number of pairs of eigenvalues of $\tilde A$ lying in the unit circle.
\end{Remark}


\subsection{A physics application}\label{phys}
Our results allow a number of interesting physics applications. Here
we mention the application of Theorem \ref{1} and Theorem
\ref{1general} to non-Hermitian crystals. While Hermiticity lies at
the heart of quantum mechanics, recent experimental advances in
controlling dissipation have brought about unprecedented flexibility
in engineering non-Hermitian Hamiltonians in open classical and
quantum systems \cite{Gong}. Non-Hermitian Hamiltonians exhibit rich
phenomena without Hermitian analogues: e.g. parity-time
($\mathcal{PT}$) symmetry breaking, topological phase transition,
non-Hermitian skin effects, e.t.c \cite{AGU,BBK}, and all of these
phenomena can be observed in  non-Hermitian crystals \cite{longhi,jiang2019interplay}.

Here, we consider the non-Hermitian crystals of the form \begin{equation}\label{non-hermitian-op}
(H_{V,x+i\e,\alpha}u)_n=u_{n+1}+u_{n-1}+V(x+i\e+n\alpha)u_n.
\end{equation}
This defines a non-self adjoint operator on $\ell^2(\Z)$. An important
class of non-Hermitian Hamiltonians which have recently attracted a significant attention in physics is called parity-time ($\mathcal{PT}$) symmetry Hamiltonian (i.e. $\overline{v}(n)=v(-n)$, \cite{BB}). Indeed, if $V$ is even with $x=0$, \eqref{non-hermitian-op} is a $\mathcal{PT}$ symmetry Hamiltonian. Different from  the self-adjoint operators,  the spectra of non-self-adjoint operators may not always consist of real numbers, and physicists  are interested in the phase transition from real energy spectrum (unbroken $\mathcal{PT}$ phase) to complex energy spectrum (broken $\mathcal{PT}$ phase), i.e. $\mathcal{PT}$ symmetry breaking  phase transition \cite{longhi}.  As first discovered in \cite{Liu2020L}, this kind of transition can be studied through the analysis of Lyapunov exponents $L_{\e}(E)$. 
Theorem \ref{1general} allows to easily  deduce that  subcritical
radius  $$ \min_{E \in\Sigma_{V,\alpha}} h(E)=
\frac{1}{2\pi} \min_{E \in\Sigma_{V,\alpha}}  \hat{L}_1(E)$$ is
the $\mathcal{PT}$ symmetry breaking parameter (one may consult
\cite{Liu2020L,LZC} for the detailed reasoning).

  Another way to understand parity-time ($\mathcal{PT}$) symmetry breaking is  topological phase transition.
  Let  $E_B\in \R$ be a base energy which is not in the  spectrum of
  $H.$ We  introduce a  \textit{ topological winding number} as
\begin{equation}
\nu(E_B,\e)= \lim_{\epsilon\rightarrow 0}\lim_{N\rightarrow\infty} \frac{1}{2\pi i}\frac{1}{N} \int_{0}^{2\pi} \partial _{\theta }\ln
\det [H_N(\theta, \e+\epsilon)-E_{B}] d\theta,
\end{equation}
where $H_N= P_{[1,N]}HP_{[1,N]}$.
The winding number $\nu$ counts the number of times the complex
spectral trajectory encircles the base point $E_B$ when the real phase
$\theta$ varies from zero to $2\pi$ \cite{Gong,longhi}.  It was shown
in \cite{LZC,WWYZ} that  \textit{topological winding number} is
precisely equal to the \textit{acceleration}:
\begin{equation} \label{windingpsi}
\nu(E_B,\e)= - \frac{1}{2\pi} \frac{\partial L_{\e}(E)}{\partial \e}.
\end{equation}
Note  that the fact that $E_B$ doesn't belongs to the  spectrum of
$H_{V,x+i\e,\alpha}$ just means that $\e$ is not a turning point of
$L_{\e}(E_B)$.


For a concrete example, one can take a non-Hermitian SEM
\begin{equation*}
(H_{\alpha,x+i\epsilon}u)(n)=u(n+1)+u(n-1)+2\lambda_1\cos2\pi(x+i\e+n\alpha)u(n)+2\lambda_2\cos2\pi(x+i\e +n\alpha)u(n).
\end{equation*}
As a consequence of \eqref{windingpsi} and Theorem \ref{1}, we have the following characterization of its  topological winding number:
 \begin{equation}
\nu (E_B, \e)= \left\{
\begin{array}{cc}
0,& ~ 0<\e< \frac{  \hat{L}_1(E_B) }{2\pi} ,\\
-1, &~ \frac{  \hat{L}_1(E_B) }{2\pi} <\e< \frac{  \hat{L}_2(E_B) }{2\pi} ,\\
-2, &~ \e>\frac{  \hat{L}_2(E_B) }{2\pi}. \\
\end{array}%
\right.
\end{equation}%
where $\hat{L}_2(E)\geq  \hat{L}_1(E) \geq 0$ are the Lyapunov exponents
of the dual operator \eqref{dual_gamo}.   One can consult \cite{LZC} for more detail.

\section{Our approach}\label{ideas}

Once discovered and formulated, the multiplicative Jensen's formula can ostensibly be
proved in different ways, some possibly being a matter of pure technique. Here,
however, we believe our method itself is almost as valuable as the
resulting formula, as we develop a dynamical perspective on the
non-self-adjoint duality, several components of which are very general
and of independent interest.

While if trying to mimic the Aubry-Andre-Avron-Kotani-Simon-Haro-Puig
approach, one can still define the IDS and prove a non-self-adjoint
Thouless formula for ergodic Schr\"odinger operators (also in the
strip) following \cite{WWYZ}, it is not clear if the invariance of the IDS holds.

Our approach starts instead with the invariance of the Green's function: 
$$
 \int_{\T}\langle \delta_0,(H_{V(\cdot+i\e),x,\alpha}-E)^{-1}\delta_0\rangle dx = \int_{\T}\langle \delta_0, (\widehat{H}_{V(\cdot+i\e),\theta,\alpha}-EI)^{-1}\delta_0\rangle d\theta
$$
Other than the Thouless formula, another important link between the
Lyapunov exponent and operator-theoretic properties of $H$ is  the Johnson-Moser's theorem \cite{JM}:
$$
\frac{\partial L_{0}(E)}{\partial \Im  E}= - \Im \int_{\T} \langle \delta_0,(H_{V,x,\alpha}-E)^{-1}\delta_0\rangle dx,
$$
which connects the derivative of the Lyapunov exponent and the
averaged Green's function. The big advantage is that it has a non-self-adjoint version,
\begin{equation*}
\frac{\partial L_{\e}(E)}{\partial \Im  E}= - \Im \int_{\T} \langle \delta_0,(H_{V(\cdot+i\e),x,\alpha}-E)^{-1}\delta_0\rangle dx
\end{equation*}
and also the strip version for individual distinct Lyapunov exponents (counting multiplicity),
\begin{equation}\label{indi}
2\pi \frac{\partial(\sum_{j=n_{i-1}+1}^{n_i}\gamma_{j})}{\partial \Im E}(E)=\frac{-1}{d}\tr\Im\int_{\T}G_{i}(\theta,E) d\theta.
\end{equation}
Finally, we develop a new general method to calculate the Green's
function of strip operators in a purely dynamical way.  This enables us to link the dual averaged Green's function
$$
 \int_{\T}\langle \delta_0, (\widehat{H}_{V(\cdot+i\e),\theta,\alpha}-EI)^{-1}\delta_0\rangle d\theta
$$
to the sums of individual averaged Green's functions in \eqref{indi},
which then links the derivative of $L_\e(E)$ and the derivative of the
right hand side of \eqref{gne1}.

Overall, our approach has three  key ingredients, each of independent
value and the last two also of a significantly higher generality
\begin{enumerate}
\item {\it Partial hyperbolicity of the dual cocycle}
  (Corollaries \ref{dominate1},  \ref{m'}).  It
  turns out that duals of Schr\"odinger cocycles, are
  either trivial or hyperbolic or partially hyperbolic, and in fact, a
   stronger domination statement holds (Theorem \ref{dom}). 
  Note  that dynamics of partially hyperbolic diffeomorphisms with 1D
  (or 2D)-center, is an important and  difficult topic in ergodic
  theory and smooth dynamical systems \cite{ACW,ACW2,AV,RH,RHRHU}.
This crucial discovery here in particular confirms the importance of
the acceleration, which is exactly half the dimension of the center, on the dynamical systems side,  and  is also
important to our further results on the Cantor spectrum \cite{martini} and
sharp phase transition conjecture \cite{gj} for type I operators. We expect it to play a central role in investigating other spectral problems.
\item {\it  Johnson-Moser's theorem for Schr\"odinger operators on the
    strip} (Proposition \ref{mg}). We develop a purely dynamical
  method to prove the classical Johnson-Moser's theorem. This method
  is of high generality and works for any finite-range operator whose cocycle is partially hyperbolic.
Our method gives the relation between the individual Lyapunov exponents
and the Green function, - a correspondence which was not known
before. This has already allowed the first author \cite{ge2} to solve a major
open problem formulated by Kotani and Simon \cite{ks} on partial reflectionless
of the M matrices of strip operators in presence  of some positive
Lyapunov exponents \footnote{We would like to thank Professor Kotani
  who pointed this out  to us.}.
\item {\it A representation of the Green’s function for general strip
    operators} (Lemma \ref{basic}). We develop a way to construct the
  Green's function of the strip operator via the half-line decaying
  solutions in a pure dynamical  way. The key is that our
  method effectively works for any {\it non-self-adjoint} operator. For example, we apply it to construct the Green's function for the complexified Schr\"odinger operators and its dual strip operators which are out of reach via spectral methods.
\end{enumerate}

\section{Preliminaries}\label{pre}
\subsection{Complex one-frequency cocycles}\label{com}
 Let $(\Omega,\tilde{d})$ be a compact metric space with metric
 $\tilde{d}$, $T:\Omega\rightarrow \Omega$ a  homeomorphism,  and let ${\rm M}_m(\C)$ be the set of all $m\times m$ matrices. Given any  $A\in C^0(\Omega,{\rm M}_m(\C))$, a cocycle $(T, A)$ is a linear skew product:
$$
(T,A)\colon \left\{
\begin{array}{rcl}
\Omega \times \C^{m} &\to& \Omega \times \C^{m}\\[1mm]
(x,v) &\mapsto& (Tx,A(x)\cdot v)
\end{array}
\right.  .
$$
For $n\in\mathbb{Z}$, $A_n$ is defined by $(T,A)^n=(T^n,A_n).$ Thus $A_{0}(x)=id$,
\begin{equation*}
A_{n}(x)=\prod_{j=n-1}^{0}A(T^{j}x)=A(T^{n-1}x)\cdots A(Tx)A(x),\ for\ n\ge1,
\end{equation*}
and $A_{-n}(x)=A_{n}(T^{-n}x)^{-1}$.

Here we are mainly interested in the case where $\Omega=\T$ is the torus, and $Tx=x+\alpha$, where $\alpha \in \R\backslash
\Q$ is an irrational  number. We call  $(\alpha,A)$ a
 {\it complex one-frequency cocycle}.
We denote by $L_1(\alpha, A)\geq L_2(\alpha,A)\geq...\geq L_m(\alpha,A)$ the Lyapunov exponents of $(\alpha,A)$ repeated according to their multiplicities, i.e.,
$$
L_k(\alpha,A)=\lim\limits_{n\rightarrow\infty}\frac{1}{n}\int_{\T}\ln(\sigma_k(A_n(x)))dx,
$$
where for any matrix $B\in {\rm M}_m(\C)$, we denote by
$\sigma_1(B)\geq...\geq \sigma_m(B)$ its singular values (eigenvalues
of $\sqrt{B^*B}$).  Note that since the k-th exterior product
$\Lambda^kB$ of $B$ satisfies $\sigma_1(\Lambda^kB)=\|\Lambda^kB\|$,
we have that $L^k(\alpha, A)=\sum\limits_{j=1}^kL_j(\alpha,A)$ satisfies
\begin{equation}\label{LEs}
L^k(\alpha,A)=\lim\limits_{n\rightarrow \infty}\frac{1}{n}\int_{\T}\ln(\|\Lambda^kA_n(x)\|)dx.
\end{equation}
Note that for $A\in C^0(\T,{\rm GL}_m(\C)),$
where ${\rm GL}_m(\C)$ is the set of all $m\times m$ invertible
matrices, we have $L_m(\alpha,A)>-\infty.$
\begin{Remark} We note that the order we choose here, as well as in
  the proofs in the next two sections is $L_1(\alpha, A)\geq
  L_2(\alpha,A)\geq...\geq L_m(\alpha,A)$ while we use the
  opposite order when we talk about dual Lyapunov exponents in the
  context of Theorem \ref{1general}.
  \end{Remark}

A basic fact about  {\it complex one-frequency cocycles} is  continuity
of the  Lyapunov exponents:
\begin{Theorem}[\cite{ajs,bj,jks}]\label{lyacon}
The functions $\R \times C^{\omega}(\T, {\rm M}_m(\C))\ni
(\alpha,A)\mapsto L_k(\alpha,A)\in [-\infty,\infty)$
are continuous at any $(\alpha',A')$ with $\alpha'\in \R\backslash\Q$.
\end{Theorem}

\begin{Remark}\label{sym}
If $A\in {\rm Sp_{2d}}(\C)$ where ${\rm Sp_{2d}}(\C)$ denotes the set of  $2d\times 2d$ complex symplectic matrices, then the Lyapunov exponents come in pairs $\{\pm L_i(\alpha,A)\}_{i=1}^d$.
\end{Remark}

\subsection{Uniform hyperbolicity and dominated splitting}\label{secds}
Given any  $A\in C^0(\Omega,{\rm Sp_{2d}}(\C))$,  we say the  cocycle $(T, A)$ is {\it uniformly hyperbolic} if for every $x \in \Omega$, there exists a continuous splitting $\C^{2m}=E^s(x)\oplus E^u(x)$ such that for some constants $C>0,c>0$, and for every $n\geqslant 0$,
$$
\begin{aligned}
\lvert A_n(x)v\rvert \leqslant Ce^{-cn}\lvert v\rvert, \quad & v\in E^s(x),\\
\lvert A_n(x)^{-1}v\rvert \leqslant Ce^{-cn}\lvert v\rvert,  \quad & v\in E^u(T^nx).
\end{aligned}
$$
This splitting is invariant by the dynamics, which means that for every $x \in \Omega$, $A(x)E^{\ast}(x)=E^{\ast}(Tx)$, for $\ast=s,u$.   The set of uniformly hyperbolic cocycles is open in the $C^0$-topology.

A related concept, 
{\it dominated splitting}, is defined the following way. Recall that
for complex one-frequency cocycles $(\alpha,A)\in C^0(\T,{\rm
  M}_m(\C))$  Oseledets theorem provides us with  strictly decreasing
sequence of Lyapunov exponents $L_j \in [-\infty,\infty)$ of multiplicity $m_j\in\N$, $1\leq j \leq \ell$ with $\sum_{j}m_j=m$, and for $a.e.$ $x$, there exists
a measurable invariant decomposition $$\C^m=E_x^1\oplus E_x^2\oplus\cdots\oplus E_x^\ell$$ with $\dim E_x^j=m_j$ for $1\leq j\leq \ell$ such that $$
\lim\limits_{n\rightarrow\infty}\frac{1}{n}\ln\|A_n(x)v\|=L_j,\ \  \forall v\in E_x^j\backslash\{0\}.
$$
An invariant decomposition $\C^m=E_x^1\oplus E_x^2\oplus\cdots\oplus E_x^\ell$  is  {\it dominated} if for any unit vector $v_j\in E_x^j\backslash \{0\}$, we have $$\|A_n(x)v_j\|>\|A_n(x)v_{j+1}\|.$$
Oseledets decomposition is a priori only measurable, however if an invariant decomposition $\C^m=E_x^1\oplus E_x^2\oplus\cdots\oplus E_x^\ell$  is  {\it dominated}, then $E_x^j$ depends  continuously on $x$ \cite{bdv}.

A cocycle $(\alpha,A)$     is called $k$-dominated (for some $1\leq k\leq m-1$) if there exists a dominated decomposition $\C^m=E^+ \oplus E^- $    with $\dim E^+ = k.$  If  $\alpha\in\R\backslash\Q$, then it follows from the definitions that the Oseledets splitting is dominated if and only if $(\alpha,A)$ is $k$-dominated for each $k$ such that  $L_k(\alpha,A)> L_{k+1}(\alpha,A)$.

\subsection{Global theory of one-frequency quasiperiodic cocycles}

We briefly introduce the global  theory of one-frequency quasiperiodic cocycles, first developed for $SL(2,\C)$-cocycles \cite{avila}, and
later  generalized to any ${\rm M}_m(\C)$-cocycles \cite{ajs}.  The most  important  concept of the global theory is the acceleration.  If $A\in C^{\omega}(\T,{\rm M}_m(\C))$ admits a holomorphic extension to $|\Im z|<\delta$, then for $|\e|<\delta$ we can define $A_\e\in C^{\omega}(\T,M_m(\C))$ by $A_\e(x)=A(x+i\e)$.
The accelerations of $(\alpha,A)$  are defined as
$$
\omega^k(\alpha,A)=\lim\limits_{\e\rightarrow 0^+}\frac{1}{2\pi\e}(L^k(\alpha,A_\e)-L^k(\alpha,A)), \qquad  \omega_k(\alpha,A)= \omega^k(\alpha,A)-\omega^{k-1}(\alpha,A).
$$
The key ingredient to the global theory is that the acceleration  is quantized.

\begin{Theorem}[\cite{avila,ajs}]\label{ace}
There exist $1\leq l\leq m$, $l\in\N$, such that $l\omega^k$ and $l \omega_k$ are integers. In particular, if $A\in C^\omega(\T, SL(2,\C))$, then $\omega^1(\alpha,A)$ is an integer.
\end{Theorem}
\begin{Remark}\label{r1}
If $L_j(\alpha,A)>L_{j+1}(\alpha,A)$, then $\omega^j(\alpha,A)$   is
an integer. This is contained in the proof of Theorem 1.4 in \cite{ajs}, see also footnote 17 in \cite{ajs}.
\end{Remark}

By subharmonicity, we know that $L^k(\alpha,A(\cdot+i\e))$ is a convex function of $\e $  in a neighborhood of $0$,  unless it is identically equal to $-\infty$.
We say that $(\alpha,A)$ is $k$-regular if $\e\rightarrow L^k(\alpha,A(\cdot+i\e))$ is an affine function of $\e$ in a neighborhood of $0$.  In general, one can relate regularity and dominated splitting as follows.

\begin{Theorem}[\cite{avila,ajs}]\label{t2.1}
Let $\alpha\in\R\backslash\Q$ and $A\in C^\omega(\T,M_m(\C))$. If $1\leq j\leq m-1$ is such that $L_j(\alpha,A)>L_{j+1}(\alpha,A)$, then $(\alpha,A)$ is $j$-regular if and only if $(\alpha,A)$ is $j$-dominated. In particular, if  $A\in C^\omega(\T, SL(2,\C))$ with $L(\alpha,A)>0$, then $(\alpha,A)$ is $1$-regular (or regular) if and only if $(\alpha,A)$ is uniformly hyperbolic.
\end{Theorem}

\subsection{Schr\"odinger operators and  Schr\"odinger cocycles}

Let $(\Omega,\tilde{d})$ be a compact metric space with distance $\tilde{d}$, $T:\Omega\rightarrow \Omega$ a  homeomorphism, and $V:\Omega\rightarrow \C$ a complex-valued continuous function. $(\Omega,T)$ is said to be {\it minimal} if each $T$-orbit is dense. We consider the following complex-valued dynamically defined Schr\"odinger operators:
\begin{equation}
(H_{x}u)_n=u_{n+1}+u_{n-1}+V(T^nx)u_n,\ \ n\in\Z,
\end{equation}
and denote by $\Sigma_x$  the spectrum of $H_{x}$. We have the following:
\begin{Lemma}\label{notde}
There is some $\Sigma\subset \C$ such that $\Sigma_x=\Sigma$ for all $x\in\Omega$.
\end{Lemma}
\begin{Remark}
  This is a standard fact for real-valued $V$ (so self-adjoint
  $H$). We provide here a brief argument that does not require self-adjointness.
\end{Remark}
\begin{pf}
We only need to prove that for any $x,y\in\Omega$, $\Sigma_x=\Sigma_y$. Assume $E\notin \Sigma_x$, that is $(H_x-E)^{-1}$ exists and is bounded. Since $(\Omega,T)$ is minimal, there is a subsequence $\{n_i\}_{i=1}^\infty$ such that $T^{n_i}y\rightarrow x$. Since $\Omega$ is compact, $T$ is uniformly continuous which implies that $H_{T^{n_i}y}\rightarrow H_x$ in operator norm. Hence $(H_{T^{n_i}y}-E)^{-1}$ exists and is bounded for $i$ sufficiently large, which implies $E\notin \Sigma_{T^{n_i}y}$ for $i$ sufficiently large. Since $H_y$ and $H_{T^{n_i}y}$ are unitary equivalent, we have $E\notin \Sigma_{y}$, thus $\Sigma_y\subset\Sigma_x$. Similarly, $\Sigma_x\subset\Sigma_y$.
\end{pf}

Note that any formal solution $u=(u_n)_{n \in \Z}$ of $H_{x}u=E u$ satisfies
$$
\begin{pmatrix}
u_{n+1}\\
u_n
\end{pmatrix}
= A_E(T^n x) \begin{pmatrix}
u_{n}\\
u_{n-1}
\end{pmatrix},\quad \forall \  n \in \Z,
$$
 where
$$
A_E(x):=
\begin{pmatrix}
E-V(x) & -1\\
1 & 0
\end{pmatrix},   \quad E\in\R.
$$ We call $(T,A_E)$  \textit{Schr\"{o}dinger cocycles}.  The
spectrum $\Sigma$ is closely related to the dynamical behavior of  the
Schr\"{o}dinger cocycle  $(T,A_E)$.  In the self-adjoint case,
i.e. the potential $V$ is real valued, then by the celebrated  Johnson's theorem \cite{johnson},  $E\notin \Sigma$ if and only if $(T,A_E)$ is uniformly hyperbolic.
It turns out that it is not difficult to extend  Johnson's theorem
\cite{johnson} to the non self-adjoint case. We will give a proof in
the appendix.

\begin{Theorem}\label{jog} Suppose that   $V:\Omega\rightarrow \C$ a complex-valued continuous function and $(\Omega,T)$ is  minimal, then $E\notin \Sigma$ if and only if $(T,A_E)$ is uniformly hyperbolic.
\end{Theorem}

In this paper,  we are mainly interested in  the following complex-valued quasiperiodic  Schr\"odinger operators
$$
(H_{V,x+i\e,\alpha}u)_n=u_{n+1}+u_{n-1}+V(x+i\e+n\alpha)u_n,
$$
and corresponding \textit{Schr\"{o}dinger cocycles}
$(\alpha,A_E(\cdot+i\e))$. Throughout the paper, we will denote $L_\e(E)=L(\alpha,A_E(\cdot+i\e))$  for short.

\subsection{Quasiperiodic Schr\"odinger operators on the strip}\label{strip}

We recall that quasiperiodic finite-range operator
\begin{equation*}
(\widehat{H}_{V,\theta,\alpha}u)(n)=\sum\limits_{k=-d}^{d} V_k u_{n+k}+2\cos 2\pi (\theta+n\alpha)u_n, \ \ n\in\Z.
\end{equation*}
naturally induces a quasiperiodic cocycle $(\alpha,\widehat{A}_{E})$ where
$$
\widehat{A}_{E}(\theta)=\frac{1}{V_d}
\begin{pmatrix}
\begin{smallmatrix}
-V_{d-1}&\cdots&-V_1&E-2\cos2\pi(\theta)-V_0&-V_{-1}&\cdots&-V_{-d+1}&-V_{-d}\\
V_d& \\
& &  \\
& & & \\
\\
\\
& & &\ddots&\\
\\
\\
& & & & \\
& & & & & \\
& & & & & &V_{d}&
\end{smallmatrix}
\end{pmatrix}.
$$
We can write it as a second order $2d$-dimensional difference equation by introducing the auxiliary variables
$$
\vec{u}_k=(u_{dk+d-1}\ \ \cdots\ \ u_{dk+1}\ \ u_{dk})^T
$$
for $k\in \Z$. It is easy to check that $(\vec{u}_k)_k$ satisfies
\begin{equation}\label{ef}
C\vec{u}_{k+1}+B(\theta_{dk})\vec{u}_k+C^*\vec{u}_{k-1}=E\vec{u}_k,
\end{equation}
where
$$
C=\begin{pmatrix}
V_d&\cdots&V_1\\
0&\ddots&\vdots\\
0&0&V_d
\end{pmatrix},
$$
$C^*$ is the transposed and conjugated matrix of $C$,  and $B(\theta)$ is the Hermitian matrix
\begin{equation}\label{strip1}
B(\theta)=\begin{pmatrix}
2\cos 2\pi (\theta_{d-1})&V_{-1}&\cdots&V_{-d+1}\\
V_1&\ddots&\ddots&\vdots\\
\vdots&\ddots&2\cos 2 \pi (\theta_1)&V_{-1}\\
V_{d-1}&\cdots&V_1&2\cos2 \pi(\theta)
\end{pmatrix}
\end{equation}
where $\theta_j=\theta+j\alpha$. Note that equation \eqref{ef} is an eigenequation of the following Schr\"odinger operator on the strip
$$
(H_{C,B,\theta,d\alpha}\vec{u})_k=C\vec{u}_{k+1}+B(T^k\theta)\vec{u}_k+C^*\vec{u}_{k-1},
$$
acting on $\ell^2(\Z,\C^d)$, which is an ergodic operator with the
dynamics given by $T\theta=\theta+d\alpha$.

To obtain a first order system and the corresponding linear skew product we use the fact that $C$ is invertible (since $V_d\neq 0$ because the degree of $V$ is exactly $d$) and write
$$
\begin{pmatrix}
\vec{u}_{k+1}\\
\vec{u}_k
\end{pmatrix}
=\begin{pmatrix}
C^{-1}(EI-B(T^k\theta))& -C^{-1}C^*\\
I_d&O_d
\end{pmatrix}
\begin{pmatrix}
\vec{u}_k\\
\vec{u}_{k-1}
\end{pmatrix}
$$
where $I_d$ and $O_d$ are the d-dimensional identity and zero matrices, respectively.
Denote
\begin{equation}\label{strip2}
\widehat{A}_{d,E}(\theta)=\begin{pmatrix}
C^{-1}(EI-B(\theta))& -C^{-1}C^*\\
I_d&O_d
\end{pmatrix}
\end{equation}
An important ingredient for our results is the complex symplectic structure
$$
\Omega=\begin{pmatrix}
0&-C^*\\
C&0
\end{pmatrix}
$$
which satisfies $\Omega^*=-\Omega$, and the fact that our Schr\"odinger skew-products $(d\alpha,\widehat{A}_{d,E})$ are complex symplectic for real $E$ with respect to $\Omega$.   However, if $E$ is complex, then $ (d\alpha,\widehat{A}_{d,E})$ are not  complex symplectic anymore.
For more details, see \cite{haro}.

We denote
$\gamma_{i}(E)=\frac{L_i(\alpha,\widehat{A}_{E})}{2\pi}$  for
$1\leq i\leq d$ for short. Then using the Kotani-Simon \cite{ks}
version of the  Thouless formula for the strip, one can prove the
following beautiful Haro-Puig's formula 
\begin{Theorem}[\cite{haro}]\label{thou}
For any $E\in \C$, we have
\be\label{lya-eq}
L(E)=2\pi\left(\sum\limits_{i=1}^{d} \gamma_i(E)\right)+\ln |V_{d}|.
\ee
\end{Theorem}

\section{Green's function of finite-range  Schr\"odinger operator}\label{hermit}
In this section, we explore the M matrix and Green's matrix for the following Schr\"odinger operators on the strip
\begin{equation}\label{strip1}
(H_{C,B,\theta,d\alpha}\vec{u})_k=C\vec{u}_{k+1}+B(T^k\theta)\vec{u}_k+C^*\vec{u}_{k-1},
\end{equation}
where $T^k\theta=\theta+kd\alpha$.
\subsection{The M matrix and the Green's matrix}

In the following, $\C_+$ will denote $\{E\in\C|\Im E>0\}$. For any energy in $\C_+$, one can define the M matrix and Green's matrix with the help of the following result:

\begin{Lemma}[\cite{haro},\cite{ks}]\label{m}
For any $E\in\C_+$, there exist unique sequences of $d\times d$ matrix
valued functions $\{F_{\pm}(k,\theta,E)\}_{k\in\Z}$ that
satisfy the following properties:
\begin{enumerate}
\item $F_\pm(0,\theta,E)=I_d,$
\item  $$
C^*F_{\pm}(k-1,\theta,E)+CF_\pm(k+1,\theta,E)+B(T^k\theta)F_\pm(k,\theta,E)=EF_\pm(k,\theta,E),
$$
\item
$$
\sum\limits_{k=0}^\infty\|F_+(k,\theta,E)\|^2<\infty, \ \ \sum\limits_{k=-\infty}^{0}\|F_-(k,\theta,E)\|^2<\infty.
$$
\end{enumerate}
\end{Lemma}

Once we have $F_{\pm}(k,\theta,E)$, we can define the M matrices
$$
M_+(\theta,E)=-F_+(1,\theta,E),
$$
$$
M_-(\theta,E)=-F_-(-1,\theta,E),
$$
as in  \cite{ks}, and note that the $M$ matrices satisfy  the following Ricatti equations.
\begin{Lemma}For any $n\in\Z$, we have
\begin{align}\label{rce}
CM_+(T^{n}\theta,E)+C^*M_+^{-1}(T^{n-1}\theta,E)+(E-B(T^n\theta))=0.
\end{align}
\begin{align}\label{rce2}
C^*M_-(T^n\theta,E)+CM_-^{-1}(T^{n+1}\theta,E)+(E-B(T^n\theta))=0.
\end{align}
\end{Lemma}
\begin{pf}
We only prove \eqref{rce}, since \eqref{rce2} can be proved similarly.  Note that $F_+(n,\theta,E)$ satisfies
$$
C^*F_+(n-1,\theta,E)+CF_+(n+1,\theta,E)+B(T^n\theta)F_+(n,\theta,E)=EF_+(n,\theta,E),
$$
and by the fact that
$$ F_+(m+n,\theta,E)= F_+(m, T^n \theta,E)F_+(n,\theta,E),$$
we have
$$
C^*M^{-1}_+(T^{n-1}\theta,E)+CM_+(T^n\theta,E)+(E-B(T^n\theta))=0.
$$
\end{pf}

Similarly as  in  \cite{ks}, one can define the Green's Matrix as
$$G(\theta,E)= \langle \vec{\delta}_0,(H_{C,B,\theta,d\alpha}-E)^{-1}\vec{\delta}_0\rangle,
$$
where
$$
\vec{\delta}_j(n)=
\begin{cases}
0& n\neq j\\
I_d &n=j
\end{cases}.
$$
The Green's matrix can be expressed as the following:

\begin{Lemma}\label{Green_Matrix}
For any $E\in\C_+$, we have
\begin{align*}
G(\theta,E)=(-CM_+(\theta,E)-C^*M_-(\theta,E)+B(\theta)-E)^{-1}
\end{align*}
\end{Lemma}
\begin{pf}
It is easy to check that
\begin{align*}
&\langle \vec{\delta}_m,(H_{C,B,\theta,d\alpha}-E)^{-1}\vec{\delta}_n\rangle \\
=&\begin{cases}
F_+(m,\theta,E)(CF_+(n+1,\theta,E)+C^*F_-(n-1,\theta,E)+(B(T^{n}\theta)-E)F_+(n,\theta,E))^{-1}& m\geq n\\
F_-(m,\theta,E)(CF_+(n+1,\theta,E)+C^*F_-(n-1,\theta,E)+(B(T^{n}\theta)-E)F_+(n,\theta,E))^{-1}& m< n
\end{cases}.
\end{align*}
%
\end{pf}

The following Lemma gives the relation between the M matrix and the Green's matrix.
\begin{Lemma} For any $E\in\C_+$, the following relation holds:
\begin{eqnarray}
\nonumber G(\theta,E)&=&(-CM_+(\theta,E)+CM^{-1}_-(T\theta,E))^{-1}, \\
\nonumber G(T\theta,E)&=& (C^*M_+^{-1}(\theta,E)-C^*M_-(T\theta,E))^{-1},\\
\label{g3} G(\theta,E)CM^{-1}_-(T\theta,E) &=& M_-(T\theta,E)G(T\theta,E) C^* + I_d.
\end{eqnarray}
\end{Lemma}

\begin{pf}
By Lemma \ref{Green_Matrix}, \eqref{rce} and \eqref{rce2}, one has
\begin{align*}
G(T\theta,E)&=(-CM_+(T\theta,E)-C^*M_-(T\theta,E)+B(T\theta)-E)^{-1}\\ \nonumber
&=(C^*M_+^{-1}(\theta,E)-C^*M_-(T\theta,E))^{-1}.
\end{align*}
\begin{align*}
G(\theta,E)&=(-CM_+(\theta,E)-C^*M_-(\theta,E)+B(\theta)-E)^{-1}\\ \nonumber
&=(-CM_+(\theta,E)+CM^{-1}_-(T\theta,E))^{-1}.
\end{align*}
Consequently, we have the following
\begin{align*}
G(\theta,E)CM^{-1}_-(T\theta,E)
=&(I_d - M_-(T\theta,E)M_+(\theta,E))^{-1}\\ \nonumber
=&M_+^{-1}(\theta,E)(M_+^{-1}(\theta,E)- M_-(T\theta,E))^{-1}\\ \nonumber
=&M_-(T\theta,E)(M^{-1}_+(\theta,E)-M_-(T\theta,E))^{-1} + I_d\\ \nonumber
=& M_-(T\theta,E)G(T\theta,E) C^* + I_d.
\end{align*}

\end{pf}

\subsection{A dynamical consequence: dominated splitting}\label{dominate}
\noindent

For any $E\in \C$, we group $\{\gamma_i(E)\}_{i=1}^d$ as
$\gamma_{n_1}, \cdots, \gamma_{n_\ell}$ with multiplicities
$\{n_{i}-n_{i-1}\}_{i=1}^\ell$ respectively, where $n_0=0$ and we assume that
$$
\gamma_{n_1}>\gamma_{n_2}>\cdots>\gamma_{n_{\ell}}\geq 0.
$$
Note that $(d\alpha,\widehat{A}_{d,E})$ is the $d$-th iteration of
$(\alpha,\widehat{A}_{E}),$ we have $L_i(d\alpha,\widehat{A}_{d,E})=2\pi d\gamma_i(E)$. Hence $\{L_{n_i}(d\alpha,\widehat{A}_{d,E})\}_{i=1}^\ell$ are the distinct Lyapunov exponents of $(d\alpha,\widehat{A}_{d,E})$. Note that we always have
$$n_{\ell(E)}(E)=d.$$
Indeed, if $E\in \R$, $(d\alpha,\widehat{A}_{d,E})$ are complex symplectic, thus
by Remark \ref{sym},  the Lyapunov exponents of $(\alpha,\widehat{A}_{E})$ come into pairs $\{\pm \gamma_i(E)\}_{i=1}^d$.  If
 $E\in\C \backslash \R$, $(d\alpha,\widehat{A}_{d,E})$ is uniformly hyperbolic \cite{haro}, thus  $L_d(d\alpha,\widehat{A}_{d,E})>0$.
 A key observation of our proof is the following:

\begin{Proposition}\label{dominate}
For $E\in\C_+$, the cocycle $(d\alpha,\widehat{A}_{d,E})$ is $n_i$-dominated for $1\leq i\leq \ell$.
\end{Proposition}
\begin{Remark}By Theorem \ref{t2.1}, this is essentially the same as part 1 of Theorem \ref{dom}.\end{Remark}
\begin{pf}
 We divide the proof into two cases:

\textbf{Case 1: $i=\ell$}. Since $E\in\C_+$, we have $(d\alpha,\widehat{A}_{d,E})$ is uniformly hyperbolic \cite{haro}, which implies that $n_{\ell}=d$ and  $(d\alpha,\widehat{A}_{d,E})$ is $d$-dominated.

\textbf{Case 2: $1\leq i\leq \ell-1$}.  Recall that
\begin{equation*}
\widehat{A}_{d,E}(\theta)=\begin{pmatrix}
C^{-1}(EI-B(\theta))& -C^{-1}C^*\\
I_d&O_d
\end{pmatrix}
\end{equation*}
where
\begin{equation*}
B(\theta)=\begin{pmatrix}
2\cos(\theta_{d-1})&V_{-1}&\cdots&V_{-d+1}\\
V_1&\ddots&\ddots&\vdots\\
\vdots&\ddots&2\cos(\theta_1)&V_{-1}\\
V_{d-1}&\cdots&V_1&2\cos(\theta)
\end{pmatrix}.
\end{equation*}
Thus if we let $\left(\ell_{ij}\right)_{1\leq i,j\leq
  d}=(\widehat{A}_{d,E})_n(\theta)$, then it is easy
to check that each $\ell_{ij}$ is a polynomial of $\cos2\pi(\theta)$
with degree $\leq n$. Similarly, if we let $L_{ij}$ be the $ij$-th entry of $\Lambda^{n_i}(\widehat{A}_{d,E})_n(\theta)$, by the definition of wedge, $L_{ij}$ is a polynomial of $\cos2\pi(\theta)$ with degree $\leq nn_i$. Hence one can compute
\begin{align*}
&|\omega^{n_i}(d\alpha,\widehat{A}_{d,E})|=\left|\lim\limits_{\e\rightarrow 0^+}\frac{1}{2\pi\e}(L^{n_i}(d\alpha,\widehat{A}_{d,E}(\cdot+i\e))-L^{n_i}(d\alpha,\widehat{A}_{d,E})\right|\\
&=\frac{1}{2\pi}\left|\lim\limits_{n\rightarrow \infty}\frac{1}{n}\int_{\T}\ln(\|\Lambda^{n_i}(\widehat{A}_{d,E})_n(\theta+i\e)\|)d\theta-\lim\limits_{n\rightarrow \infty}\frac{1}{n}\int_{\T}\ln(\|\Lambda^{n_i}(\widehat{A}_{d,E})_n(\theta)\|)d\theta\right|
\leq n_i.
\end{align*}
It follows that
$$
|\omega^{n_i}(\alpha,\widehat{A}_{E})|=\left|\frac{\omega^{n_i}(d\alpha,\widehat{A}_{d,E})}{d}\right|\leq \frac{n_i}{d}<1.
$$

On the other hand, since $\gamma_{n_i}(E)>\gamma_{n_{i}+1}(E)$, by
Remark  \ref{r1}, $\omega^{n_i}(\alpha,\widehat{A}_{E})$ is an
integer. Together with the fact that $|\omega^{n _i}(\alpha,\widehat{A}_{E})|$ is strictly smaller than $1$, we have $\omega^{n_i}(\alpha,\widehat{A}_{E})=0$.  This implies that
$$
L^{n_i}(\alpha,\widehat{A}_{E}(\cdot+i\e))=L^{n_i}(\alpha,\widehat{A}_{E})
$$
for $\e>0$ which is sufficiently small. Similar  argument  works for  $\e<0$  which is also sufficiently small.  This means  $(\alpha,\widehat{A}_{E})$ is $n_i$-regular.
Notice that
$$
dL^{n_i}(\alpha,\widehat{A}_{E}(\cdot+i\e))=L^{n_i}(d\alpha,\widehat{A}_{d,E}(\cdot+i\e)).
$$
so $(\alpha,\widehat{A}_{E})$ is $n_i$-regular if and only if  $(d\alpha,\widehat{A}_{d,E})$ is $n_i$-regular. Hence
 $(d\alpha,\widehat{A}_{d,E})$ is $n_i$-dominated by Theorem \ref{t2.1}.
\end{pf}

\begin{Corollary}\label{dominate1}
Assume $\alpha\in \R\backslash\Q,$ $E\in\R$, and $\gamma_d=0.$ Then
$(d\alpha,\widehat{A}_{d,E})$ is partially hyperbolic with
center of dimension $2(n_\ell-n_{\ell-1})$.
\end{Corollary}
\begin{pf}
For any $E\in \R$, by the same argument as in Proposition \ref{dominate}, $(d\alpha,\widehat{A}_{d,E})$ is $n_i$-dominated for any $1\leq i\leq \ell-1$. Together with the fact that $(d\alpha,\widehat{A}_{d,E})$ is complex symplectic, we have $(d\alpha,\widehat{A}_{d,E})$ is partially hyperbolic with a center of dimension $2(n_\ell-n_{\ell-1})$.
\end{pf}
\begin{Remark} ByTheorem \ref{t2.1}, this also provides a proof of part 2 of Theorem \ref{dom}.\end{Remark}
{\bf Proof of Corollary \ref{ph}:} It follows directly from the
combination of  Corollaries \ref{3general} and \ref{dominate1}. \qed

\smallskip
Notice that Proposition \ref{dominate} gives an enhancement of Lemma \ref{m}. As a consequence, we obtain:
\begin{Corollary}\label{m'}
For any $E\in\C_+$, there are sequences of $d\times d$ matrix valued functions $\{\widetilde{F}_{\pm}(k,\theta,E)\}_{k\in\Z}$  and $B(\theta,E)\in C^\omega(\T\times \C_+, \rm{GL}_{d\times d}(\C))$,  with $
\widetilde{F}_\pm(0,\theta,E)=B(\theta,E),
$  obeying 
$$
C^*\widetilde{F}_{\pm}(k-1,\theta,E)+C\widetilde{F}_\pm(k+1,\theta,E)+B(T^k\theta)\widetilde{F}_\pm(k,\theta,E)=E\widetilde{F}_\pm(k,\theta,E),
$$
$$
\sum\limits_{k=0}^\infty\|\widetilde{F}_+(k,\theta,E)\|^2<\infty, \ \ \sum\limits_{k=-\infty}^{0}\|\widetilde{F}_-(k,\theta,E)\|^2<\infty.
$$
Moreover, if we denote
\begin{align*}
\widetilde{F}_\pm(k,\theta,E)
=\begin{pmatrix}\vec{f}^\pm_1(k,\theta,E)&\vec{f}^\pm_2(k,\theta,E)&\cdots&\vec{f}^\pm_{d}(k,\theta,E)\end{pmatrix},
\end{align*}
then for any $\theta\in\T$,
$$
\limsup\limits_{k\rightarrow \infty}\frac{1}{2k}\ln\left(\|\vec{f}^-_j(k,\theta,E)\|^2+\|\vec{f}^-_j(k+1,\theta,E)\|^2\right)=2\pi d\gamma_{n_i}(E),\ \  n_{i-1}+1\leq j\leq n_i,
$$
\end{Corollary}

\begin{pf}By Proposition \ref{dominate},   $(d\alpha,\widehat{A}_{d,E})$ is $n_i$-dominated for $1\leq i\leq \ell$.  Thus by    Appendix B in   \cite{bdv},  there exist continuous invariant decompositions $E_s(\theta)$ and  $E_u(\theta)=E_u^1(\theta)\oplus E_u^2(\theta)\oplus\cdots\oplus E_u^\ell(\theta)$ such that $\C^{2d}=E_s(\theta)\oplus E_u(\theta)$.  Moreover, we have
\begin{equation}\label{eq10'}
\limsup\limits_{k\rightarrow\infty}\frac{1}{k}\ln\|(\widehat{A}_{d,E})_k(\theta)\vec{v}\|=d\gamma_{n_i}(E)>0,\ \  \forall \vec{v}\in E_u^i(\theta)\backslash\{0\},\ \   \forall \theta\in\T.
\end{equation}

Note  that $E_s(\theta)$ and $\{E^i_u(\theta)\}_{i=1}^\ell$ depend continuously on $\theta$. Actually, by Theorem 6.1 in \cite{ajs}, if the cocycle is analytic, $E_s(\theta)$ and $\{E^i_u(\theta)\}_{i=1}^\ell$ can be chosen to depend holomorphically on both $E\in \C_+$ and $\theta\in\T$. i.e., there exists $\widetilde{F}_-(\theta,E)\in C^\omega(\T\times \C_+, \rm{M}_{2d\times d}(\C))$ and $M_-^i(\theta,E)\in C^\omega(\T\times\C_+,\rm{GL}_{n_i-n_{i-1}}(\C))$ such that
\begin{align}\label{M1}
(\widehat{A}_{d,E}(\theta))^{-1}\widetilde{F}_-(\theta,E)=\widetilde{F}_-(T^{-1}\theta,E)\text{diag} \left\{-M_-^1(\theta,E),-M_-^2(\theta,E),\cdots, -M_-^\ell(\theta,E)\right\}.
\end{align}
Now, we define
$$
\begin{pmatrix}\widetilde{F}_-(k,\theta,E)\\
\widetilde{F}_-(k-1,\theta,E)\end{pmatrix}=(\widehat{A}_{d,E})_k(\theta)\widetilde{F}_-(\theta,E).
$$
It follows from \eqref{eq10'}, for any $\theta\in\T$,
\begin{equation}\label{new1}
\limsup\limits_{k\rightarrow \infty}\frac{1}{2k}\ln\left(\|\vec{f}^-_j(k,\theta,E)\|^2+\|\vec{f}^-_j(k+1,\theta,E)\|^2\right)=2\pi d\gamma_{n_i}(E),\ \  n_{i-1}+1\leq j\leq n_i,
\end{equation}
Finally, we take $B(\theta,E)=\widetilde{F}_-(0,\theta,E)$ and
$$
 \widetilde{F}_+(\theta,E)=\begin{pmatrix}I_d\\
F_+(-1,\theta,E)\end{pmatrix}\cdot B(\theta,E), \qquad
\begin{pmatrix}\widetilde{F}_+(k,\theta,E)\\
\widetilde{F}_+(k-1,\theta,E)\end{pmatrix}=(\widehat{A}_{d,E})_k(\theta)\widetilde{F}_+(\theta,E).
$$
Then by \eqref{new1} and Lemma \ref{m}, we have
$$
\sum\limits_{k=0}^\infty\|\widetilde{F}_+(k,\theta,E)\|^2<\infty, \ \ \sum\limits_{k=-\infty}^{0}\|\widetilde{F}_-(k,\theta,E)\|^2<\infty.
$$
\end{pf}

\begin{Remark}\label{relation}
By the uniqueness of $F_\pm(k,\theta,E)$, it is standard that $\widetilde{F}_\pm(k,\theta,E)$ and $F_\pm(k,\theta,E)$ have the following relations
\begin{equation*}
\widetilde{F}_\pm(k,\theta,E)=F_\pm(k,\theta,E)B(\theta,E).
\end{equation*}
\end{Remark}

It follows directly that
\begin{Corollary}\label{Green_Matrix'}
For any $E\in\C_+$, we have
\begin{align*}
G(\theta,E)=\widetilde{F}_-(0,\theta,E)(C\widetilde{F}_+(1,\theta,E)-C\widetilde{F}_-(1,\theta,E))^{-1}.
\end{align*}
\end{Corollary}
\begin{pf}
By Lemma \ref{Green_Matrix} and Remark \ref{relation}, we have
\begin{align*}
G(\theta,E)&=F_+(0,\theta,E)(CF_+(1,\theta,E)+C^*F_-(-1,\theta,E)+(B(\theta)-E)F_+(0,\theta,E))^{-1}\\
&=F_+(0,\theta,E)(CF_+(1,\theta,E)-CF_-(1,\theta,E))^{-1}\\
&=\widetilde{F}_+(0,\theta,E)B^{-1}(\theta,E)(C\widetilde{F}_+(1,\theta,E)B^{-1}(\theta,E)-C\widetilde{F}_-(1,\theta,E)B^{-1}(\theta,E))^{-1}\\
&=\widetilde{F}_-(0,\theta,E)(C\widetilde{F}_+(1,\theta,E)-C\widetilde{F}_-(1,\theta,E))^{-1}.
\end{align*}
\end{pf}

Corollary \ref{m'} also implies that  the $M_-(\theta,E)$ is conjugated to a block diagonal matrix.
\begin{Proposition}\label{f1}
There exist $B\in C^\omega(\T\times\C_+,{\rm GL_d}(\C))$ and $M_-^i(\theta,E)\in C^\omega(\T\times\C_+,{\rm GL_{n_i-n_{i-1}}}(\C))$ for $1\leq i\leq \ell$ such that
$$
\widetilde{M}_-(\theta,E):=B^{-1}(T^{-1}\theta,E)M_-(\theta,E)B(\theta,E)=\text{diag} \{M_-^1(\theta,E),M_-^2(\theta,E),\cdots, M_-^\ell(\theta,E)\}.
$$
Moreover, for $1\leq i\leq \ell$, if we denote by
$$
\omega_i(E)=\int_\T\ln\det{M_-^i(\theta,E)}d\theta,
$$
then
\be\label{reo}
\Re \omega_i(E)=2\pi d(n_{i}-n_{i-1})\gamma_{n_{i}}(E).
\ee
\end{Proposition}

\begin{pf}Notice  that \eqref{M1} and Remark \ref{relation} imply
\begin{align}\label{M3}
\nonumber \begin{pmatrix}F_-(1,\theta,E)\\ F_-(0,\theta,E) \end{pmatrix} B(\theta,E)=&\begin{pmatrix}F_-(0,T^{-1}\theta,E)\\ F_-(-1,T^{-1}\theta,E) \end{pmatrix} B(T^{-1}\theta,E)\\
& \text{diag} \left\{-M_-^1(\theta,E),-M_-^2(\theta,E),\cdots, -M_-^\ell(\theta,E)\right\}.
\end{align}
It follows that
\begin{align}\label{M3}
\nonumber \begin{pmatrix}F_-(0,\theta,E)\\ F_-(-1,\theta,E) \end{pmatrix}B(\theta,E)=&\begin{pmatrix}F_-(1,T^{-1}\theta,E)\\ F_-(0,T^{-1}\theta,E) \end{pmatrix} B(T^{-1}\theta,E)\\
&\text{diag} \left\{-M_-^1(\theta,E),-M_-^2(\theta,E),\cdots, -M_-^\ell(\theta,E)\right\}.
\end{align}
Thus we have
$$
B^{-1}( T^{-1}\theta,E)M_-(\theta,E)B(\theta,E)=\text{diag} \{M_-^1(\theta,E),M_-^2(\theta,E),\cdots, M_-^\ell(\theta,E)\}.
$$
For \eqref{reo}, we only prove the case $i=1$, the others follow similarly. Note
\begin{align*}
&2\pi dn_{1}\gamma_{n_{1}}(E)=\lim\limits_{n\rightarrow \infty}\frac{1}{n}\int_{\T}\ln\left(\|\Lambda^{n_1}(\widehat{A}_{d,E})_n(\theta)\vec{f_1}(0,\theta,E)\wedge\cdots\wedge\vec{f}_{n_1}(0,\theta,E) \|\right)d\theta\\
=&\lim\limits_{n\rightarrow \infty}\frac{1}{n}\int_{\T}\ln\left(\|\vec{f_1}(n-1,\theta,E)\wedge\cdots\wedge\vec{f}_{n_1}(n-1,\theta,E) \|\right)d\theta\\
=&\lim\limits_{n\rightarrow \infty}\frac{1}{n}\left(\sum_{j=0}^{n-1}\ln|\det{M_{-}^1(\theta+j\alpha,E)}|+\int_{\T}\ln\left(\|\vec{f_1}(0,T^{n-1}\theta,E)\wedge\cdots\wedge\vec{f}_{n_1}(0,T^{n-1}\theta,E) \|\right)d\theta\right)\\
=&\int_\T\ln|\det{M_-^1(\theta,E)}|d\theta= \Re \omega_1(E).
\end{align*}
\end{pf}

For any $E\in \C\backslash \R,$ the classical Thouless formula will
imply Johnson-Moser's type result: $$\frac{\partial
  L^d(\alpha,\widehat{A}_{E})}{\partial \Im E}=-\frac{1}{d}\Im
\tr\int  G(\theta,E)d\theta.$$ We refer readers to \cite{haro,ks} for
details. We now provide the following generalized version of the
Thouless formula for a Lyapunov-invariant subspace.  Denote by $P_I$
the projection from $\Z$ to $I.$ We have the following generalization of Johnson-Moser's theorem:
\begin{Proposition}\label{mg}
For $1\leq i\leq \ell$, we have
$$ \frac{\partial\omega_i}{\partial E}(E)=2\pi \frac{\partial(\sum_{j=n_{i-1}+1}^{n_i}\gamma_{j})}{\partial \Im E}(E)=\frac{-1}{d}\tr\Im\int_{\T}G_{i}(\theta,E) d\theta.$$
where $G_i(\theta)=P_{[n_{i-1}+1,n_{i}]}B^{-1}(\theta,E)G(\theta,E)B(\theta,E)P_{[n_{i-1}+1,n_{i}]}$.
\end{Proposition}

We postpone the proof of  Proposition \ref{mg} to Section \ref{proofmg}.

\section{Green's function for non-self-adjoint quasiperiodic  operators}\label{non-hermit}
We start with establishing Aubry duality between a non-Hermitian quasiperiodic Schr\"odinger operator on $\ell^2(\Z)$:
\begin{equation}\label{csch}
(H_{V(\cdot+i\e),x,\alpha}u)(n)=u_{n+1}+u_{n-1}+V(x+i\e+n\alpha)u_n, \ \ n\in\Z,
\end{equation}
and  the finite range quasiperiodic  operator  $\hat{L}_{V(\cdot+i\e),\theta,\alpha}$:
\begin{equation}\label{strip1}
(\widehat{H}_{V(\cdot+i\e),\theta,\alpha}u)(n)=\sum\limits_{k=-d}^{d} e^{-k\e}V_k u_{n+k}+2\cos2\pi(\theta+n\alpha)u_n, \ \ n\in\Z,
\end{equation}
where $V(x)=\sum\limits_{k=-d}^{d} V_ke^{2\pi ikx}$ is a trigonometric
polynomial.  Then we will analyze the Green's function for these non-Hermitian quasiperiodic operators.

\subsection{Aubry duality}

Consider the fiber direct integral,
$$
\mathcal{H}:=\int_{\T}^{\bigoplus}\ell^2(\Z)dx,
$$
which, as usual, is defined as the space of $\ell^2(\Z)$-valued, $L^2$-functions over the measure space $(\T,dx)$.  The extensions of the
Sch\"odinger operators  and their long-range duals to  $\mathcal{H}$ are given in terms of their direct integrals, which we now define.
Let $\alpha\in\T$ be fixed. Interpreting $H_{V(\cdot+i\e),x,\alpha}$ as fibers of the decomposable operator,
$$
H_{V(\cdot+i\e),\alpha}:=\int_{\T}^{\bigoplus}H_{V(\cdot+i\e),x,\alpha}dx,
$$
the family $\{H_{V(\cdot+i\e),x,\alpha}\}_{x\in\T}$ naturally induces an operator on the space $\mathcal{H}$, i.e. ,
$$
(H_{V(\cdot+i\e),\alpha} \Psi)(x,n)= \Psi(x,n+1)+ \Psi(x,n-1) +  V(x+i\e+n\alpha) \Psi(x,n).
$$

Similarly,  the direct integral of finite-range operator  $\widehat{H}_{V(\cdot+i\e),\theta,\alpha}$,
denoted as $\widehat{H}_{V(\cdot+i\e),\alpha}$, is given by
$$
(\widehat{H}_{V(\cdot+i\e),\alpha}  \Psi)(x,n)=  \sum\limits_{k=-d}^{d} e^{-2\pi k\e}V_k \Psi(x,n+k)+  2\cos2\pi (\theta+n\alpha) \Psi(x,n).
$$
These two operators are bounded and  non-Hermitian in $\mathcal{H}$.  Let us now see that  operators \eqref{csch} and \eqref{strip1} are in fact unitarily equivalent.
Indeed, by analogy with the heuristic and classical approach to Aubry duality \cite{gjls}, let  $U$ be the following operator on $\mathcal{H}:$
$$
(\mathcal{U}\phi)(\eta,m):=\hat{\phi}(m, \eta+\pi\alpha m)=\sum_{n\in\Z}\int_{\T}e^{2\pi imx}e^{2\pi in(m\alpha+\eta)}\phi(x,n)dx.
$$
$U$ is clearly unitary and a direct computation show that it
conjugates $H$ and $\hat{L}$
$$U H_{V(\cdot+i\e),\alpha} U^{-1}=  \widehat{H}_{V(\cdot+i\e),\alpha}.$$
Moreover, we have the following:
\begin{Lemma}\label{12}
For any $z\in\C\backslash\R$, one has
$$
\int_{\T}\langle (H_{V(\cdot+i\e),x,\alpha}-z)^{-1}\delta_n,\delta_n\rangle dx=\int_{\T}\langle (\widehat{H}_{V(\cdot+i\e),\theta,\alpha}-z)^{-1}\delta_n,\delta_n\rangle d\theta.
$$
\end{Lemma}
\begin{pf}
Let $\phi(\theta,m)=\delta_n$ for any $\theta\in\T.$ By the unitary equivalence between  $H_{V(\cdot+i\e),\alpha}$ and $\hat{L}_{V(\cdot+i\e),\alpha}$, we have that
\begin{align*}
\int_{\T}\langle (H_{V(\cdot+i\e),\theta,\alpha}-z)^{-1}\delta_n,\delta_n\rangle d\theta&=\int_{\T}\langle (H_{V(\cdot+i\e),\alpha}-z)^{-1} \phi(\theta,m),\phi(\theta,m)\rangle d\theta\\
&=\int_{\T}\langle \mathcal{U}(H_{V(\cdot+i\e),\alpha}-z)^{-1} \phi(\theta,m),\mathcal{U}\phi(\theta,m)\rangle d\theta\\
&=\int_{\T}\langle \mathcal{U} (H_{V(\cdot+i\e),\alpha}-z)^{-1} \mathcal{U}^{-1}\mathcal{U}\phi(\theta,m),\mathcal{U}\phi(\theta,m)\rangle d\theta\\
&=\int_{\T}\langle (\widehat{H}_{V(\cdot+i\e),\theta,\alpha}-z)^{-1}\delta_0,\delta_0\rangle d\theta\\
&=\int_{\T}\langle (\widehat{H}_{V(\cdot+i\e),\theta,\alpha}-z)^{-1} \delta_n,\delta_n\rangle d\theta
\end{align*}
where we used the fact $(\mathcal{U}\phi)(\theta,m)=e^{2\pi in\theta}
\delta_0$.
\end{pf}

\subsection{A representation formula for the Green's function}
In this subsection, we state a general lemma which is useful for the computation of Green's function of  finite-range operators. Its proof can be found in Section \ref{proofbasic}.
\begin{Lemma}\label{basic}
Consider the following $2d$ order difference operator,
$$
(Lu)(n)=\sum\limits_{k=-d}^da_ku(n+k)+V(n)u(n).
$$
If the eigenequation $Lu=Eu$ has $2d$ linearly independent solutions $\{\phi_i\}_{i=1}^{2d}$ satisfying
$$
\phi_i\in\ell^2(\Z^+)(i=1,\cdots,m)\ \ \phi_i\in\ell^2(\Z^-)(i=m+1,\cdots,2d),
$$
then $L-EI$ is invertible. Moreover,
$$
\langle\delta_p,(L-EI)^{-1}\delta_q\rangle=\begin{cases}
\frac{\sum\limits_{i=1}^m\phi_i(p)\Phi_{1,i}(q)}{a_d\det{\Phi(q)}} &\text{$p\geq q+1$}\\
-\frac{\sum\limits_{i=m+1}^{2d}\phi_i(p)\Phi_{1,i}(q)}{a_d\det{\Phi(q)}} &\text{$p\leq q$}
\end{cases},
$$
where
$$
\Phi(q)=\begin{pmatrix}\phi_1(q+d)&\phi_2(q+d)&\cdots&\phi_{2d}(q+d)\\ \phi_1(q+d-1)&\phi_2(q+d-1)&\cdots&\phi_{2d}(q+d-1)\\ \vdots&\vdots& &\vdots\\
\phi_1(q-d+1)&\phi_2(q-d+1)&\cdots&\phi_{2d}(q-d+1)\end{pmatrix}
$$
and $\Phi_{i,j}(q)$ is the $(i,j)$-th cofactor of $\Phi(q)$.
\end{Lemma}

\begin{Remark}\label{2d+}
If for the eigenequation $Lu=Eu$ there exist $2d$ independent solutions $\{\phi_i\}_{i=1}^{2d}$ with
$
\phi_i\in\ell^2(\Z^+)(i=1,\cdots,2d)$,
then  we have
$$
\langle\delta_p,(L-EI)^{-1}\delta_q\rangle=\begin{cases}
\frac{\sum\limits_{i=1}^m\phi_i(p)\Phi_{1,i}(q)}{a_d\det{\Phi(q)}} &\text{$p\geq q+1$}\\
0 &\text{$p\leq q$}
\end{cases}.
$$

\end{Remark}

\subsection{Green's function of non-Hermitian Schr\"odinger operators}

In this subsection, we study the Green function of \eqref{csch} for $\e \neq 0$. In this case, we first have the following:

\begin{Lemma}\label{muse}
If $(\alpha,A_E(\cdot+i\e))$ is regular and $L_\e(E)>0$, then there are unique solutions $u_\pm(k,x+i\e,E)\in\ell^2(\Z^\pm)$, obeying
$$
u_\pm(k-1,x+i\e,E)+u_\pm(k+1,x+i\e,E) + V(x+i\e+k\alpha)u_\pm(k,x+i\e,E)=Eu_\pm(k,x+i\e,E),
$$
where $
u_\pm(0,x+i\e,E)=1.
$
\end{Lemma}
\begin{pf}
By Theorem \ref{t2.1}, $(\alpha, A_E(\cdot +i\e))$ is uniformly hyperbolic. The existence of $u_\pm$ follows from the definition of uniform hyperbolicity.
\end{pf}

Once we have this, similar to the Hermitian case, we can define the  $m$ function
as
$$
m_{\pm}(x+i\e,E)= - u_{\pm}(\pm1,x+i\e,E),
$$
and one can express the Green's function  defined as
\begin{align*}
g(x+i\e,E)=\langle \delta_0,(H_{V(\cdot+i\e),x,\alpha}-E)^{-1}\delta_0\rangle,
\end{align*}
by the m-function as follows:

\begin{Lemma}\label{non-sch}
$g(x+i\e,E)=\frac{-1}{m_+(x+i\e,E)+m_-(x+i\e,E)+E-V(x+i\e)}.$
\end{Lemma}
\begin{pf}
By Proposition \ref{basic} and Lemma \ref{muse}, we have
\begin{align*}
g(x+i\e,E)&=\frac{1}{u_+(1,x+i\e,E)-u_-(1,x+i\e,E)}\\
&=\frac{1}{u_+(1,x+i\e,E)+u_-(-1,x+i\e,E)+V(x+i\e)-E}\\
&=\frac{-1}{m_+(x+i\e,E)+m_-(x+i\e,E)+ E-V(x+i\e)}.
\end{align*}
\end{pf}

One can now relate the derivative of Lyapunov exponent and Green function as follows:

\begin{Proposition}\label{sch}
If $(\alpha,A_E(\cdot+i\e))$ is regular and $L_\e(E)>0$, then
$$
\frac{\partial L_{\e}(E)}{\partial \Im  E}= - \Im \int_{\T} g(x+i\e,E) dx.
$$
\end{Proposition}
\begin{pf}
Similarly to the Hermitian case, $m$-function is non zero for any
$x\in\T$, so we can define
  $$w_\e(E)=\int \ln m_-(x+i\e,E) dx.$$
Thus it suffices for us to prove
\be\label{schle}
\frac{\partial w_{\e}(E)}{\partial E}= \int_{\T} g(x+i\e,E) dx.
\ee
Once we have this, then the result follows from the Cauchy-Riemann equation directly.

To prove \eqref{schle}, first note that we also have the Ricatti equation
\begin{eqnarray}
\nonumber && m_+(x+\alpha+i\e,E)+m_+^{-1}(x+i\e,E)+(E-V(x+\alpha+ i\e))=0,\\
\label{mp}
&&m_-(x+i\e,E)+m_-^{-1}(x+\alpha+i\e,E)+(E-V(x+i\e))=0.
\end{eqnarray}
By Lemma \ref{non-sch}, this implies that
\begin{align}\label{ft-2}
\begin{split}
g(x+i\e,E)&= \frac{-1}{m_+(x+i\e,E)-m^{-1}_-(x+\alpha+i\e,E)},\\
g(x+\alpha+i\e,E) &= \frac{-1}{m_-(x+\alpha+i\e,E)-m^{-1}_+(x+i\e,E)},\\
g(x+i\e,E)m_+(x+i\e,E) &=
g(x+\alpha+i\e)m_-(x+\alpha+i\e,E).
\end{split}
\end{align}

 We now introduce the   auxiliary function \footnote{A similar idea will be used again in the proof of  Proposition \ref{mg}}
$$f(x,E)=g(x+i\e,E)\frac{\partial m_-(x+i\e,E)}{ \partial E}.$$

By taking the derivative of \eqref{mp}, we have
$$
\frac{\partial m_-(x+i\e,E)}{ \partial E}=\frac{1}{m_-^{2}(x+\alpha+i\e,E)}\frac{\partial m_-(x+\alpha+i\e,E)}{ \partial E}-1.
$$
Then by \eqref{ft-2}, it follows that
\begin{align*}
&f(x+\alpha,E)-f(x,E)\\
=&g(x+\alpha+i\e,E)\frac{\partial m_-(x+\alpha+i\e,E)}{ \partial E}-g(x+i\e,E)\left(\frac{1}{m_-^{2}(x+\alpha+i\e,E)}\frac{\partial m_-(x+\alpha+i\e,E)}{ \partial E}-1\right)\\
=&\left(g(x+\alpha+i\e,E)-g(x+i\e,E)\frac{1}{m_-^{2}(x+\alpha+i\e,E)}\right)\frac{\partial m_-(x+\alpha+i\e,E)}{ \partial E}+g(x+i\e,E)\\
=& g(x+\alpha+i\e,E)\left(1- \frac{1}{m_+(x+i\e,E)m_-(x+\alpha+i\e,E)}\right) \frac{\partial m_-(x+\alpha+i\e,E)}{ \partial E}+g(x+i\e,E)\\
=& -\frac{1}{m_-(x+\alpha+i\e,E)}\frac{\partial m_-(x+\alpha+i\e,E)}{ \partial E}+g(x+i\e,E).
\end{align*}
Taking the integral over $\T$, we get the desired result.
\end{pf}

\subsection{Green's function for the duals of non-Hermitian Schr\"odinger operators}

In this subsection, we study  the Green's function of  the operator \eqref{strip1}.
Note $\widehat{H}_{V(\cdot+i\e),\theta,\alpha}$ naturally induces a quasiperiodic cocycle  $(\alpha,\widehat{A}^{\e}_{E})$ where
$$
\widehat{A}^{\e}_{E}(x)=\frac{1}{e^{-d\e}V_d}
\begin{pmatrix}
\begin{smallmatrix}
-e^{-2\pi(d-1)\e}V_{d-1}&\cdots&-e^{-2\pi\e}V_1&E-2\cos2\pi(x)-V_0&-e^{2\pi\e} V_{-1}&\cdots&-e^{2\pi(d-1)\e}V_{-d+1}&-e^{2\pi d\e}V_{-d}\\
e^{-2\pi d\e}V_d& \\
& &  \\
& & & \\
\\
\\
& & &\ddots&\\
\\
\\
& & & & \\
& & & & & \\
& & & & & &e^{-2\pi d\e}V_{d}&
\end{smallmatrix}
\end{pmatrix}.
$$
For $1\leq i\leq 2d,$ we denote by $\gamma_{i}^{\e}(E)=L_i(\alpha,\widehat{A}^{\e}_{E})$  for short. The basic observation is the following:

\begin{Lemma}\label{9}
We have
\begin{equation}\label{new}
\gamma^\e_i(E)=\gamma_i(E)+\e , \quad 1\leq i\leq 2d.
\end{equation}
\end{Lemma}

\begin{pf} By a direct computation, one can prove that
\begin{equation} \label{rn1}
\widehat{A}^{\e}_{E}=e^{2\pi \e} D_{d}^{-1}\widehat{A}_{E}D_{d},
\end{equation}
where
$$D_d =diag\{e^{2\pi d\e}, e^{2\pi(d-1)\e},\cdots, e^{-2\pi(d-2)\e}, e^{-2\pi(d-1)\e}\}. $$
Thus \eqref{new} follows from  the definition of Lyapunov exponents.
\end{pf}

With this observation in hand, one can express the Green's  function of \eqref{strip1}.  Indeed, for any $E\in\C_+$, recall that we assume the Lyapunov exponent of $(\alpha,
\widehat{A}_{E})$ satisfy
\be \label{lya-re}
\gamma_{n_1}>\gamma_{n_2}>\cdots>\gamma_{n_{\ell}}>0
\ee
with multiplicity of each $\gamma_{n_i}$ being $\{n_{i}-n_{i-1}\}_{i=1}^\ell$.   For simplicity of the notations, in the following,    we will  just denote $\gamma_{n_0} =\infty$, $\gamma_{n_{\ell+1}}=0$, and rewrite \eqref{lya-re} as
$$
 \infty =\gamma_{n_0} >\gamma_{n_1}>\gamma_{n_2}>\cdots>\gamma_{n_{\ell}}>  \gamma_{n_{\ell+1}}=0.
$$

Thus we can give the following representation of Green's function from Proposition \ref{basic}.
\begin{Proposition}\label{green function}
For any fixed $E\in\C_+$,  $0\leq i\leq \ell$,   if $\e\in (-\gamma_{n_{i}}(E),-\gamma_{n_{i+1}}(E))$, then $\widehat{H}_{V(\cdot+i\e),\theta,\alpha}-EI$ is invertible for any $\theta\in\T$. Moreover, we have
$$
\int_{\T}\langle \delta_0, (\widehat{H}^{2\cos}_{V(\cdot+i\e),\theta,\alpha}-EI)^{-1}\delta_0\rangle d\theta =\begin{cases}
 \frac{1}{d}\sum\limits_{j=1}^{i}\int_{\T}\tr G_j(\theta,E)d\theta  &\text{$1 \leq i\leq \ell$}\\
0 &\text{$i=0$}
\end{cases}.
$$
\end{Proposition}
\begin{pf}

First note that $u_n$ solves \begin{equation*}
(\widehat{H}_{V,\theta,\alpha}u)(n)=\sum\limits_{k=-d}^{d} V_k u_{n+k}+2\cos2\pi(\theta+n\alpha)u_n, \ \ n\in\Z
\end{equation*}
if and only if
$$
\vec{u}_k=(u_{dk+d-1}\ \ \cdots\ \ u_{dk+1}\ \ u_{dk})^T
$$  solves
\begin{equation*}
C\vec{u}_{k+1}+B(T^k\theta)\vec{u}_k+C^*\vec{u}_{k-1}=E\vec{u}_k.
\end{equation*}
Thus by Corollary \ref{m'}, $\widetilde{F}_\pm(k,\theta,E)$ can be written as
\begin{align*}
\widetilde{F}_\pm(k,\theta,E)&=\begin{pmatrix}\vec{f}^\pm_1(k,\theta,E)&\vec{f}^\pm_2(k,\theta,E)&\cdots&\vec{f}^\pm_{d}(k,\theta,E)\end{pmatrix}\\
&=\begin{pmatrix}f^\pm_1(kd+d-1,\theta,E)&f^\pm_2(kd+d-1,\theta,E)&\cdots&f^\pm_{d}(kd+d-1,\theta,E)\\ f^\pm_1(kd+d-2,\theta,E)&f^\pm_2(kd+d-2,\theta,E)&\cdots&f^\pm_{d}(kd+d-2,\theta,E)\\ \vdots&\vdots& &\vdots\\
f^\pm_1(kd,\theta,E)&f^\pm_2(kd,\theta,E)&\cdots&f^\pm_{d}(kd,\theta,E)\end{pmatrix} \nonumber
\end{align*}
and $\{f_j^\pm(n,\theta,E)\}_{j=1}^d$ are $2d$ linearly independent solutions of $\widehat{H}_{V,\theta,a}u=Eu$. Furthermore,  we have
$$
\limsup\limits_{n\rightarrow \infty}\frac{1}{2n}\ln\left(\|\vec{f}^-_j(n+1,\theta,E)\|^2+\|\vec{f}^-_j(n,\theta,E)\|^2\right)=2\pi d\gamma_{n_i}(E),\ \  n_{i-1}+1\leq j\leq n_i,
$$
$$
f_j^+(n,\theta,E)\in\ell^2(\Z^+), \ \ 1\leq j\leq d.
$$
By \eqref{rn1}, it is obvious that $\{e^{n\e}f_j^\pm(n,\theta,E)\}_{j=1}^d$ are $2d$ independent solutions of $\widehat{H}_{V(\cdot+i\e),\theta,\alpha}u=Eu$. Thus for $\e\in (-\gamma_{n_{i}}(E),-\gamma_{n_{i+1}}(E))$ and for any $\theta\in\T$, we have
\begin{align}\label{so1}
e^{2\pi n\e}f_j^+(n,\theta,E)\in\ell^2(\Z^+), \ \ 1\leq j\leq d,
\end{align}
\begin{align}\label{so2}
e^{2\pi n\e}f_j^-(n,\theta,E)\in\ell^2(\Z^-), \ \ 1\leq j\leq n_i,
\end{align}
\begin{align}\label{so3}
e^{2\pi n\e}f_j^-(n,\theta,E)\in\ell^2(\Z^+), \ \ n_{i}+1\leq j\leq d.
\end{align}
We divide the proof into two cases:\\

\textbf{Case I: $i=0$}. In this case, $e^{n\e}f_j^+(n,\theta,E)\in\ell^2(\Z^+)$, $e^{n\e}f_j^-(n,\theta,E)\in\ell^2(\Z^+)$, where $1\leq j\leq d$. Then the result follows directly from Proposition \ref{basic} (see also Remark \ref{2d+}). \\

\textbf{Case II: $1\leq i\leq \ell$.}
In this case, we first denote
$$
\Phi(n,\theta,E)=
\begin{pmatrix}
\begin{smallmatrix}
f_1^+(n+d,\theta,E)&\cdots&f_d^+(n+d,\theta,E)&f_1^-(n+d,\theta,E)&\cdots&f_d^-(n+d,\theta,E)\\ f_1^+(n+d-1,\theta,E)&\cdots&f_d^+(n+d-1,\theta,E)&f_1^-(n+d-1,\theta,E)&\cdots&f_d^-(n+d-1,\theta,E)\\ \vdots&& \vdots&\vdots&&\vdots\\
f_1^+(n-d+1,\theta,E)&\cdots&f_d^+(n-d+1,\theta,E)&f_1^-(n-d+1,\theta,E)&\cdots&f_d^-(n-d+1,\theta,E)\end{smallmatrix}
\end{pmatrix}.
$$
Let $\Phi_{i,j}(n,\theta,E)$ be the $(i,j)$-th cofactor of $\Phi(n,\theta,E)$. Then the fundamental matrix of
$$\widehat{H}_{V(\cdot+i\e),\theta,\alpha}u=Eu$$
can be rewritten as
\be \label{phie}\Phi^{\e}(n,\theta,E)=\text{diag}\{e^{2\pi(n+d)\e}, e^{2\pi(n+d-1)\e},\cdots, e^{2\pi(n-d+1)\e}\} \Phi(n,\theta,E).\ee
Let
$\Phi_{i,j}^\e(n,\theta,E)$ be the $(i,j)$-th cofactor of
$\Phi^{\e}(n,\theta,E)$. A direct computation shows that
\be\label{phiee}\Phi_{1, j}^\e(n,\theta,E) =e^{2\pi n(2d-1)\e} \Phi_{1,j}(n,\theta,E).\ee

Thus for any $\e\in (-\gamma_{n_{i}}(E),-\gamma_{n_{i+1}}(E))$, by \eqref{so1}-\eqref{phiee} and Proposition \ref{basic}, we have
\begin{align}\label{neede1}
&\sum\limits_{j=0}^{d-1}\langle \delta_j, (\widehat{H}_{V(\cdot+i\e),\theta,\alpha}-EI)^{-1}\delta_j\rangle\\  \nonumber
=&\sum_{j=0}^{d-1} \frac{-1}{V_de^{-2\pi d\e}\det{\Phi^\e(j,\theta,E)}}\sum\limits_{k=1}^{n_i} e^{2\pi j\e} f^-_k(j,\theta,E)\Phi_{1,d+k}^\e(j,\theta,E)\\ \nonumber
=&\sum_{j=0}^{d-1} \frac{-1}{V_d e^{4\pi jd\e}\det{\Phi(0,\theta,E)}}\sum\limits_{k=1}^{n_i} e^{2\pi j\e} f^-_k(j,\theta,E) e^{2\pi j(2d-1)\e}\Phi_{1,d+k}(j,\theta,E)\\ \nonumber
=&\frac{-1}{V_d \det{\Phi(0,\theta,E)}}\sum_{j=0}^{d-1}\sum\limits_{k=1}^{n_i} f^-_k(j,\theta,E)\Phi_{1,d+k}(j,\theta,E). \nonumber
\end{align}

We need the following equivalent representation of the Green's matrix.
\begin{Lemma}[Element version]\label{green2}
For any $E\in\C_+$, $\theta\in\T$ and $p,q\in\Z$, we have
$$
\langle\delta_p,(\widehat{H}_{V,\alpha,\theta}-EI)^{-1}\delta_q\rangle=\begin{cases}
\frac{\sum\limits_{i=1}^df^+_i(p,\theta,E)\Phi_{1,i}(q,\theta,E)}{V_d\det{\Phi(0,\theta,E)}} &\text{$p\geq q+1$}\\
-\frac{\sum\limits_{i=1}^{d}f^-_i(p,\theta,E)\Phi_{1,d+i}(q,\theta,E)}{V_d\det{\Phi(0,\theta,E)}} &\text{$p\leq q$}
\end{cases}.
$$
As a corollary, we have
\begin{align*}
G(\theta,E)=&\frac{-1}{V_d\det{\Phi(0,\theta,E)}}\begin{pmatrix}f^-_1(d-1,\theta,E)&f^-_2(d-1,\theta,E)&\cdots&f^-_{d}(d-1,\theta,E)\\ f^-_1(d-2,\theta,E)&f^-_2(d-2,\theta,E)&\cdots&f^-_{d}(d-2,\theta,E)\\ \vdots&\vdots& &\vdots\\
f^-_1(0,\theta,E)&f^-_2(0,\theta,E)&\cdots&f^-_{d}(0,\theta,E)\end{pmatrix}\\
&\cdot\begin{pmatrix}\Phi_{1,d+1}(d-1,\theta,E)&\Phi_{1,d+1}(d-2,\theta,E)&\cdots&\Phi_{1,d+1}(0,\theta,E)\\ \Phi_{1,d+2}(d-1,\theta,E)&\Phi_{1,d+2}(d-2,\theta,E)&\cdots&\Phi_{1,d+2}(0,\theta,E)\\ \vdots&\vdots& &\vdots\\ \Phi_{1,2d}(d-1,\theta,E)&\Phi_{1,2d}(d-2,\theta,E)&\cdots&\Phi_{1,2d}(0,\theta,E)\end{pmatrix}.
\end{align*}
\end{Lemma}
\begin{pf}
By uniqueness of the Green matrix, $G(\theta,E)$ can be written as
\begin{align*}
G(\theta,E)=\begin{pmatrix}\langle\delta_{d-1},(\hat{L}_{V,\theta,\alpha}-EI)^{-1}\delta_{d-1}\rangle&\cdots&\langle\delta_{d-1},(\hat{L}_{V,\theta,\alpha}-EI)^{-1}\delta_{0}\rangle\\ \langle\delta_{d-2},(\hat{L}_{V,\theta,\alpha}-EI)^{-1}\delta_{d-1}\rangle&\cdots&\langle\delta_{d-2},(\hat{L}_{V,\theta,\alpha}-EI)^{-1}\delta_{0}\rangle\\ \vdots&\vdots &\vdots\\
\langle\delta_{0},(\hat{L}_{V,\theta,\alpha}-EI)^{-1}\delta_{d-1}\rangle&\cdots&\langle\delta_{0},(\hat{L}_{V,\theta,\alpha}-EI)^{-1}\delta_{0}\rangle\end{pmatrix}.
\end{align*}
Note that $f_j^+(n,\theta,E)\in\ell^2(\Z^+)$ and
$f_j^-(n,\theta,E)\in\ell^2(\Z^-)$ for $1\leq i\leq d.$  Thus the
result follows from Proposition \ref{basic}, and the fact that $G(\theta,E)$ is symmetric.

\end{pf}

By Corollary \ref{Green_Matrix'} and Lemma \ref{green2}, we have
\begin{align}\label{neede2}
&(C\widetilde{F}_+(1,\theta,E)-C\widetilde{F}_-(1,\theta,E))^{-1}\\ \nonumber
=&\frac{-1}{V_d\det{\Phi(0,\theta,E)}}\begin{pmatrix}\begin{smallmatrix}\Phi_{1,d+1}(d-1,\theta,E)&\Phi_{1,d+1}(d-2,\theta,E)&\cdots&\Phi_{1,d+1}(0,\theta,E)\\ \Phi_{1,d+2}(d-1,\theta,E)&\Phi_{1,d+2}(d-2,\theta,E)&\cdots&\Phi_{1,d+2}(0,\theta,E)\\ \vdots&\vdots& &\vdots\\ \Phi_{1,2d}(d-1,\theta,E)&\Phi_{1,2d}(d-2,\theta,E)&\cdots&\Phi_{1,2d}(0,\theta,E)\end{smallmatrix}\end{pmatrix}.
\end{align}
By \eqref{neede2},  for any $1\leq i\leq \ell$, we have
\begin{align}\label{neede3}
&\frac{-1}{V_d\det{\Phi(0,\theta,E)}}\sum_{j=0}^{d-1}\sum\limits_{k=1}^{n_i} f^-_k(j,\theta,E)\Phi_{1,d+k}(j,\theta,E)\\ \nonumber
=&\frac{-1}{V_d\det{\Phi(0,\theta,E)}}\sum\limits_{k=1}^{n_i}\sum_{j=0}^{d-1} f^-_k(j,\theta,E)\Phi_{1,d+k}(j,\theta,E)\\ \nonumber
=&\tr P_{[1,n_i]}(C\widetilde{F}_+(1,\theta,E)-C\widetilde{F}_-(1,\theta,E))^{-1}\widetilde{F}_-(0,\theta,E)P_{[1,n_i]}\\ \nonumber
=&\tr P_{[1,n_i]}B^{-1}(\theta,E)G(\theta,E)B(\theta,E)P_{[1,n_i]}\\ \nonumber
=&\sum\limits_{j=1}^{i}\tr {G_j(\theta)}.
\end{align}
Using \eqref{neede1}, \eqref{neede3} and taking the integral over $\T$, we get the result.
\end{pf}

As a result of Aubry duality, we have the following corollary:

\begin{Corollary}\label{green function+}
For any fixed $E\in\C_+$, $0\leq i\leq \ell$, if $\e\in (-\gamma_{n_{i}}(E),-\gamma_{n_{i+1}}(E))$, then  Schr\"odinger cocycle $(\alpha,A_E(\cdot+i\e))$ is regular and $L_\e(E)>0$. Moreover,
\begin{align}\label{schrderivate}
\frac{\partial L_{\e}(E)}{\partial \Im  E}=\begin{cases}
 -\frac{1}{d}\sum\limits_{j=1}^{i} \tr \Im \int_{\T} G_j(\theta,E)d\theta& 1\leq i\leq \ell\\
0&i=0
\end{cases}.
\end{align}
\end{Corollary}

\begin{pf}
  For any $0\leq i\leq \ell$ and for any $\e\in (-\gamma_{n_{i}}(E),-\gamma_{n_{i+1}}(E))$, by Proposition \ref{green function}, $(E-\widehat{H}_{V(\cdot+i\e),\alpha,\theta})^{-1}$ exists and is bounded  for any $\theta\in\T$, thus $E\notin \Sigma(\widehat{H}_{V(\cdot+i\e),\alpha,\theta})$ \footnote{$\Sigma(H)$ denotes the spectrum of $H$.}, moreover by Lemma \ref{12},
  \be\label{equifreen}
 \int_{\T}\langle \delta_0,(H_{V(\cdot+i\e),x,\alpha}-E)^{-1}\delta_0\rangle dx = \int_{\T}\langle \delta_0, (\widehat{H}_{V(\cdot+i\e),\theta,\alpha}-EI)^{-1}\delta_0\rangle d\theta.
  \ee
Also by  Lemma \ref{notde},  we have $E\notin
\Sigma(H_{V(\cdot+i\e),\alpha,x})$ for any $x\in\T$.  By Theorem
\ref{jog}, $(\alpha,A_E(\cdot+i\e))$ is uniformly hyperbolic, so by Theorem \ref{t2.1}, $(\alpha,A_E(\cdot+i\e))$ is regular and $L_\e(E)>0$.  Then \eqref{schrderivate}
follows from \eqref{equifreen}, Proposition \ref{sch} and Proposition \ref{green function}.

\end{pf}

\section{The trigonometric polynomial case: Proof of Theorem \ref{1}}\label{main-proof}
 We assume that $V$ is a trigonometric polynomial of degree $d.$ For any $E\in\C$, recall that
$$
\gamma_{n_i(E)}(E):=\frac{L_{n_i(E)}(\alpha,\widehat{A}_{E})}{2\pi}, \ \ 1\leq i\leq \ell(E).
$$
We may assume that
$$
\gamma_{n_1(E)}(E)>\gamma_{n_2(E)}(E)>\cdots>\gamma_{n_{\ell(E)}(E)}(E)\geq 0,
$$
with multiplicities $\{n_{i}(E)-n_{i-1}(E)\}_{i=1}^{\ell(E)}$
respectively where $n_0(E)=0$ and, by an argument at the beginning of
Section \ref{dominate} we have $n_{\ell(E)}=d.$

While Theorem \ref{1} consider the energy $E\in \R$, we will derive it
from the following stronger result:
 \begin{Theorem}\label{7.1}
For any $\alpha\in\R\backslash\Q$ and $E\in\C$, we have the following:
\begin{align}\label{for1}
 L_{\e}(E)
     =  \left\{
\begin{aligned} & L(E) &\e\in
(-\gamma_{d}(E),0]\\
&L_{-\gamma_{n_{i+1}(E)}(E)}(E)-2\pi (d-n_{i}(E))(\e+\gamma_{n_{i+1}(E)}(E))
&\e\in (-\gamma_{n_{i}(E)}(E),-\gamma_{n_{i+1}(E)}(E)]\\
&L_{-\gamma_{1}(E)}(E)-2\pi d(\e+\gamma_{1}(E)) &\e\in (-\infty,-\gamma_{1}(E)]
\end{aligned}\right.
\end{align}
where $1\leq i\leq \ell(E)-1$.
\end{Theorem}

 \noindent \textbf{Proof of Theorem \ref{1}:} For any  $E\in\R$,  $L_\e(E)$ is an even function in $\e$. Let $\hat{L}_i(E)=2\pi \gamma_{d-i}(E)$. Theorem \ref{1} follows directly from Theorem \ref{7.1}. \qed \\

 \noindent \textbf{Proof of Theorem \ref{7.1}:} To prove  Theorem \ref{7.1}, we only need to prove the following:
\begin{Proposition}\label{7.1+}
For $\alpha\in\R\backslash\Q$ and $E\in\C$, there exists a sequence $E_n\in\C\backslash \R$, such that $E_n\rightarrow E$ and \eqref{for1} holds for each $E_n$.
\end{Proposition}

Once we have this,  Theorem \ref{7.1} can be obtained by the continuity
arguments as follows. We only prove the result for  $\e \in (-\gamma_{1}(E),-\gamma_{d}(E)]$, since the case $\e\in (-\gamma_{d}(E),0]\cup (-\infty,-\gamma_{1}(E)]$ follows directly from Proposition \ref{7.1+} and Theorem \ref{lyacon}.

For any fixed $E\in\C$, we fix $1\leq i\leq \ell(E)-1$ and $\e\in
(-\gamma_{n_i(E)}(E),-\gamma_{n_{i+1}(E)}(E))$.  By Theorem \ref{7.1},
there exists a sequence $E_n\in\C \backslash \R$ such that
$E_n\rightarrow E$ and \eqref{for1} holds for each $E_n.$ Based on  the Thouless formula (Theorem \ref{thou}),  we have
$$L(E_n)=2\pi\sum\limits_{j=1}^{d} \gamma_j(E_n)+\ln |V_{d}|=2\pi \sum\limits_{i=1}^{\ell(E_n)} (n_i(E_n)-n_{i-1}(E_n))   \gamma_{n_i(E_n)}(E_n)+\ln |V_{d}|,
$$
thus formula \eqref{for1} can be rewritten as
\begin{align}\label{for2}
L_{\e}(E_n)=\left\{
\begin{aligned} & L(E_n) &\e\in
(-\gamma_{d}(E_n),0],\\
&-2\pi(d-n_{i}(E_n))\e+2\pi\sum\limits_{j=1}^{n_i(E_n)}\gamma_{j}(E_n)+\ln|V_d|
&\e\in (-\gamma_{n_{i}(E_n)}(E_n),-\gamma_{n_{i+1}(E_n)}(E_n)],\\
&-2\pi d\e+\ln|V_d| &\e\in (-\infty,-\gamma_{1}(E_n)].
\end{aligned}\right.
\end{align}

Let $j(E_n)$ be such that $
-\gamma_{j}(E_n)< \e<-\gamma_{j+1}(E_n)$. Note that  by our selection,  $\gamma_{n_{i+1}(E)}(E)= \gamma_{n_{i}(E)+1}(E).$ Thus by continuity of
$\gamma_{n_{i}(E)}(E)$ and $\gamma_{n_{i}(E)+1}(E)$ (Theorem \ref{lyacon}),
 there exists some $N>0$, such that if $n>N$, then $j(E_n)=n_{i}(E)$ (independent of $E_n$). By \eqref{for2}, we have
\begin{align*}
L_{\tilde{\e}}(E_n)=-2\pi(d-n_i(E))\tilde{\e}+2\pi\sum\limits_{j=1}^{n_i(E)}\gamma_j(E_n)+\ln|V_d|, \quad \tilde{\e}\in (-\gamma_{n_{i}(E)}(E_n),-\gamma_{n_{i}(E)+1}(E_n)).
\end{align*}
First let $E_n\rightarrow E.$ By the continuity of Lyapunov exponents (Theorem \ref{lyacon}), we have
\begin{align}\label{ned1}
L_{\tilde{\e}}(E)=-2\pi(d-n_i(E))\tilde{\e}+2\pi\sum\limits_{j=1}^{n_i(E)}\gamma_j(E)+\ln|V_d|, \quad \tilde{\e}\in (-\gamma_{n_{i}(E)}(E),-\gamma_{n_{i}(E)+1}(E)).
\end{align}
On the other hand,  if we take
$$(E_n,\varepsilon_n)=(E_n, -\gamma_{n_{i}(E)+1}(E_n)+ \frac{\gamma_{n_{i}(E)}(E_n)-\gamma_{n_{i}(E)+1}(E_n)}{n} ),$$
then by Theorem \ref{lyacon}, $L_{\e}(E)$ is jointly continuous in
$(\e,E)$, so it follows that
\begin{align}\label{ned2}
L_{-\gamma_{n_{i}(E)+1}(E)}(E)=2\pi(d-n_{i}(E))\gamma_{n_{i}(E)+1}(E)+2\pi\sum\limits_{j=1}^{n_i(E)}\gamma_{j}(E)+\ln|V_d|.
\end{align}
By \eqref{ned1} and \eqref{ned2}, and the fact that $\gamma_{n_{i+1}(E)}(E)= \gamma_{n_{i}(E)+1}(E)$, we have
$$
L_\e(E)=L_{-\gamma_{n_{i+1}(E)}(E)}(E)-2\pi(d-n_{i}(E))(\e+\gamma_{n_{i+1}(E)}(E)), \quad \e\in (-\gamma_{n_{i}(E)}(E),-\gamma_{n_{i+1}(E)}(E)].
$$
This completes the proof. \qed

\subsection{Proof of  Proposition \ref{7.1+}}
 Proposition \ref{7.1+} follows from 
\begin{Proposition}\label{l6.1}
For  $\alpha\in\R\backslash\Q$ and $E\in\C$,   there exists  sequence $E_n\in\C \backslash \R$, such that $E_n\rightarrow E$ and
\begin{enumerate}
\item $\{-\gamma_{n_i(E_n)}(E_n)\}_{i=1}^{\ell(E_n)}$ are exactly the turning points of $L_\e(E_n)$ for $\e<0$;
\item  The variation of the slope at $-\gamma_{n_i(E_n)}(E_n)$ is $n_i(E_n)-n_{i-1}(E_n)$ for $1\leq i\leq \ell(E_n)$.
\end{enumerate}
\end{Proposition}
Indeed, for $V(x)=\sum_{k=-d}^dV_ke^{2\pi ikx}$ with $\overline{V_k}=V_{-k}$,  one has
$$L_{\e}(E)\leq \sup_{x\in \T} \ln \|A_{E}(x+i\e) \| \leq d \e + O(1).$$
Thus by convexity, for any $E\in\C$, the absolute value of the slope of $L_\e(E)$ as a function of $\e$ is less than or equal to $d$.
By a direct computation,  for  sufficiently large $\e$,
\begin{equation*}
		A_{E}(x+i\e)=e^{-2\pi d \e}e^{-2\pi i dx}\begin{pmatrix}
			-V_d&0\\0 &0
		\end{pmatrix}+ o(1).
	\end{equation*}
By the continuity of Lyapunov exponent (Theorem \ref{lyacon}), we have
	\begin{equation*}
	L_{\e}(E)= -2\pi d\e+\ln |V_d | +o(1),
	\end{equation*}
	thus by Theorem \ref{ace},
	\begin{equation*}
		L_{\e}(E)= -2\pi d\e+\ln |V_d |  \quad \text{as $\e \rightarrow -\infty$,}
	\end{equation*}
i.e. the slope of $L_\e(E)$ is $-d$, as $\e \rightarrow -\infty$. On the other hand, $L_\e(E)$ as a function of $\e$ is a piecewise convex affine function, $\sum_{i=1}^\ell \left(n_i(E_n)-n_{i-1}(E_n)\right)=d$ and there are no other turning points  except $\{-\gamma_{n_i(E_n)}(E_n)\}_{i=1}^{\ell(E_n)}$ when $\e<0$. For the piecewise affine function $L_\e(E_n)$, we know all the turning points $\{-\gamma_{n_i(E_n)}(E_n)\}_{i=1}^{\ell(E_n)}$ and the  variation of the slope at each turning point and the final slope when $\e<0$, thus we have the full information on $L_\e(E_n)$ when $\e<0$.\qed

\subsection{Proof  of Proposition \ref{l6.1}}

For simplicity, we only prove the result for $E\in \C_+ \cup \R$.  We define
$$
\mathcal{I}=\bigcup\limits_{E\in\C}\bigcup\limits_{i=1}^{\ell(E)}\left\{[n_{i-1}(E),n_i(E)]\right\},
$$
$$
\mathcal{Z}=\bigcup\limits_{I\in\mathcal{I}}\left\{E\in\C_+|\tr \Im \left(P_{I}B^{-1}(\theta,E)G(\theta,E)B(\theta,E)P_{I}\right)=0\right\}.
$$
Notice that  for any $1\leq i, j\leq d$,  $B^{-1}(\theta,E)G(\theta,E)B(\theta,E)\in C^\omega(\T\times  \C_+)$, and $\mathcal{I}$ has finitely many elements.
Thus $\mathcal{Z}$ has no cluster points.
Hence for any $E\in\C_+\cup \R$, there is a sequence $E_n\in\C_+$ with
$E_n\rightarrow E$, such that $E_n\notin \mathcal{Z},$
$$
\tr{\Im\int_\T G_i(\theta,E_n)d\theta}\neq 0,\ \ 1\leq i\leq \ell(E_n).
$$

\noindent \textbf{Proof  of Proposition \ref{l6.1} (1):} Proposition \ref{l6.1} (1) is implied directly by the following general fact:

\begin{Lemma}\label{turning}
For $E\notin \mathcal{Z} \subset \C_+$, 
$\{-\gamma_{n_i(E)}(E)\}_{i=1}^{\ell(E)}$ are exactly the turning points of $L_\e(E)$ for $\e<0$.
\end{Lemma}
\begin{pf}
  We first need the following observation
  \begin{Lemma}\label{regfor}
 If $(\alpha,A)$ is regular, then
      $\omega_-(\alpha,A')=\omega_+(\alpha,A')=\omega_+(\alpha,A)$ for
      all $A'$ in a small neighborhood of $A$, where

$$
\omega_\pm(\alpha,A)=\lim\limits_{\e\rightarrow 0^\pm}\frac{L(\alpha, A_\e)-L(\alpha,A)}{2\pi \e}.
$$

\end{Lemma}
\begin{pf}
Since $L_\e(\alpha,A)$ is convex as a function of $\e$, we have
$\omega_+(\alpha,A)$ is upper semi-continuous and $\omega_-(\alpha,A)$
is lower semi-continuous, and $(\alpha,A)$ is regular if and only if  $\omega_-(\alpha,A)=\omega_+(\alpha,A)$.  Note that regularity is an open condition in $\R\backslash\Q\times C^\omega(\T,SL(2,\C)).$
 This implies $\omega_-(\alpha,A')=\omega_+(\alpha,A')=\omega_+(\alpha,A)$ for all $A'$ in a small neighborhood of $A$.
\end{pf}

%
%
%

We will now prove Lemma \ref{turning} by contradiction.  Note that  Corollary \ref{green function+} implies that if
$$
\e\in (-\gamma_{n_{i}(E)}(E),-\gamma_{n_{i+1}(E)}(E))
$$
where $0\leq i\leq \ell(E)$, then $(\alpha, A_E(\cdot+i\e))$ is
regular, i.e. $\e$ is not a turning point. Thus we only need to prove
$\{-\gamma_{n_i(E)}(E)\}_{i=1}^{\ell(E)}$ are turning points   of
$L_\e(E)$. Otherwise, assume there exists $1\leq i_0\leq \ell(E)$ such
that $-\gamma_{n_{i_0}(E)}(E)$ is not a turning point, so
$(\alpha, A_{E}(\cdot-i\gamma_{n_{i_0}(E)}(E)))$ is regular. By
Lemma \ref{regfor}, there  exists an open rectangle $I\times J$
containing $(E,-\gamma_{n_{i_0}(E)}(E))$, such that for any
$(E',\e)\in I\times J$, there exists $m\in\Z$ (which only depends on
$(E,-\gamma_{n_{i_0}(E)}(E))$), such that the accelerations satisfy
$$ \omega_+(\alpha, A_{E'}(\cdot+i\e))=m,$$
Consequently by the same argument as in Proposition 5 in \cite{avila}, there exists $g\in C^\omega(I)$ such that for any $(E',\e)\in I\times J$, we have
$$
L_\e(E')=g(E')+m\e,
$$
which implies that
\begin{align*}
\frac{\partial L_\e}{\partial \Im E}(E')=\frac{\partial g}{\partial \Im E}(E')
\end{align*}
i.e.  $\frac{\partial L_\e}{\partial \Im E}(E')$ is independent of $\e\in J$.

On the other hand, for any $\e\in (-\gamma_{n_{i_0}(E)}(E),-\gamma_{n_{i_0+1}(E)}(E))$, by Corollary \ref{green function+}, one has
$$
\frac{\partial L_\e}{\partial \Im E}(E)=-\frac{1}{d}\sum\limits_{j=1}^{i_0}\tr \Im \int_{\T}G_{j}(\theta,E) d\theta,
$$
and  for any $\e\in (-\gamma_{n_{i_0-1}(E)}(E),-\gamma_{n_{i_0}(E)}(E))$, one has
$$
\frac{\partial L_\e}{\partial \Im E}(E) =\begin{cases} - \frac{1}{d}\sum\limits_{j=1}^{i_0-1}\tr \Im \int_{\T}G_{j}(\theta,E) d\theta &  \text{$2 \leq i_0\leq l(E)$} \\
0  & \text{$i_0=1$} \end{cases}.
$$
Thus $\frac{\partial L_\e}{\partial \Im E}(E)$ varies when $\e$ goes through $-\gamma_{n_{i_0}(E)}(E)$ since by our assumption $E\notin \mathcal{Z}.$ This  is a contradiction.

\end{pf}

\noindent \textbf{Proof  of Proposition \ref{l6.1} (2):} For any $E\in \C$ and  $1\leq i\leq \ell(E)$, we denote
$$
\omega_{n_i(E)}^+(E)=\lim\limits_{\e\searrow -\gamma_{n_i(E)}(E)}\frac{L_{\e}(E)-L_{-\gamma_{n_{i}(E)}(E)}(E)}{2\pi(\e+\gamma_{n_{i}(E)}(E))},
$$
$$
\omega_{n_i(E)}^-(E)=\lim\limits_{\e\nearrow -\gamma_{n_i(E)}(E)}\frac{L_{\e}(E)-L_{-\gamma_{n_{i}(E)}(E)}(E)}{2\pi(\e+\gamma_{n_{i}(E)}(E))}.
$$
We first prove a useful lemma
\begin{Lemma}\label{usel1}
Assume that $E\notin \mathcal{Z}$,  $1\leq i\leq \ell(E)$. If there exists $\delta>0$ such that
\be \label{lyadis}\gamma_{n_{i-1}(E)}(E')>\gamma_{n_{i-1}(E)+1}(E')=\cdots=\gamma_{n_{i}(E)}(E')>\gamma_{n_{i}(E)+1}(E') \ee
 for all $E'\in\C$ with $|E'-E|<\delta,$ then $$\omega_{n_i(E)}^+(E)-\omega_{n_i(E)}^-(E)=n_{i}(E)-n_{i-1}(E).$$
\end{Lemma}
\begin{pf}
Since there exists $\delta>0$ such that \eqref{lyadis} holds for all $E'\in\C$ with $|E'-E|<\delta$, then by the definition of $n_i$,  there exists $s(E')\in\Z$, such that
\be\label{equa-ly}
n_{i-1}(E)= n_{s-1}(E')     \qquad  n_{i}(E)= n_{s}(E'),
\ee
and one can rewrite \eqref{lyadis} as
$$\gamma_{n_{s-1}(E')}(E')>\gamma_{n_{s-1}(E')+1}(E')=\cdots=\gamma_{n_{s}(E')}(E')>\gamma_{n_{s}(E')+1}(E').$$
Without loss of generality, we can shrink $\delta$ and assume $E'\notin \mathcal{Z}$ since $\mathcal{Z}$ has at most finitely points.  Then by Lemma \ref{turning},
 $-\gamma_{n_s(E')}(E')$ is the only turning point for $\e\in (-\gamma_{n_{s-1}(E')}(E'),-\gamma_{n_{s}(E')+1}(E'))$.
If we assume
\begin{align*}
\omega_{n_i(E)}^+(E)-\omega_{n_i(E)}^-(E)=k_i(E),
\end{align*}
then by Lemma \ref{regfor}, for any $-\gamma_{n_{s-1}(E')}(E')<\e<-\gamma_{n_s(E')}(E')$, there is $m_i(E)\in\Z$ (does not depend on $E'$), such that
$$
L_{\e}(E')-L_{-\gamma_{n_{s}(E')}(E')}(E')=2\pi m_{i}(E)(\e+\gamma_{n_{s}(E')}(E')),
$$
and for any $-\gamma_{n_{s}(E')+1}(E')>\e'>-\gamma_{n_s(E')}(E')$, we have
$$
L_{\e'}(E')-L_{-\gamma_{n_{s}(E')}(E')}(E')=2\pi (m_{i}(E)+ k_i(E))(\e'+\gamma_{n_{s}(E')}(E')).
$$
Therefore, we have
$$ L_{\e'}(E')- L_{\e}(E')= 2\pi m_{i}(E) (\e'-\e)+2\pi k_i(E)  \gamma_{n_{s}(E')}(E').
$$
By Proposition \ref{mg}, one has
\begin{align*}
\frac{\partial L_{\e'}}{\partial \Im E}(E')-\frac{\partial L_{\e}}{\partial \Im E}(E')&=2\pi k_i(E)\frac{\partial\gamma_{n_s(E')}}{\partial \Im E}(E')=\frac{- k_i(E)}{d(n_s(E')-n_{s-1}(E'))}\tr\Im\int_{\T}G_{s}(\theta,E') d\theta.
\end{align*}

On the other hand, by Corollary \ref{green function+}, we have
\begin{align*}
&\frac{\partial L_{\e'}}{\partial \Im E}(E')-\frac{\partial L_{\e}}{\partial \Im E}(E')=-\frac{1}{d}\tr\Im \int_{\T}G_{s}(\theta,E') d\theta.
\end{align*}
Thus we obtain

$$
\left(\frac{k_i(E)}{n_s(E')-n_{s-1}(E')}-1\right)\frac{1}{d}\tr\Im \int_{\T}G_{s}(\theta,E') d\theta=0.
$$
By \eqref{equa-ly} and the selection $ E'\notin \mathcal{Z}$, we have $$k_i(E)=n_s(E')-n_{s-1}(E')= n_i(E)-n_{i-1}(E).$$
\end{pf}

Now we will prove that for any $1\leq i\leq \ell(E_n)$,
\begin{equation}\label{x}
\omega_{n_i(E_n)}^+(E_n)-\omega_{n_i(E_n)}^-(E_n)=n_{i}(E_n)-n_{i-1}(E_n).
\end{equation}
We distinguish  two different cases:\\

\textbf{Case I:} There exists $\delta>0$ such that $$\gamma_{n_{i-1}(E_n)}(E')>\gamma_{n_{i-1}(E_n)+1}(E')=\cdots=\gamma_{n_{i}(E_n)}(E')> \gamma_{n_{i}(E_n)+1}(E')$$ for all $E'\in\C$ with $|E'-E_n|<\delta$.  Then \eqref{x} follows directly from Lemma \ref{usel1}.\\

\textbf{Case II:} There exists $E_n^j\rightarrow E_n$ with $E_n^j \notin \mathcal{Z},$ $ 1\leq i\leq \ell(E_n^j)$,  such that not all of $$\gamma_{n_{i-1}(E_n)+1}(E^j_n),\cdots,\gamma_{n_{i}(E_n)}(E_n^j)$$ are equal. In this case, we need the following observation:

\begin{Lemma}\label{B2}
Let $a_i\in \C\rightarrow \R$ be continuous with $a_1(E)\geq
\cdots\geq a_n(E)$ for any $E\in\C$. Then for any $E_0\in\R$, there is
a sequence $E_j\rightarrow E_0$ such that for each $E_j$, there is $\delta_j>0$ and $0=i_0<i_1<\cdots<i_k=n$ with
$$
a_{i_{m-1}+1}(E)=\cdots=a_{i_m}(E), \ \  1\leq m\leq k,
$$
for any $|E-E_j|<\delta_j$.
\end{Lemma}
\begin{pf}
We only need to prove that for any $E_0\in\R$ and $j\in \Z$, there is an open set $U_j\subset B_{\frac{1}{j}}(E_0)$ and $0=i_0<i_1<\cdots<i_k=n$ such that
$$
a_{i_{m-1}+1}(E)=\cdots=a_{i_m}(E), \ \  1\leq m\leq k,
$$
for any $E\in U_j$. Then one may take $\delta_j$ such that
$B_{\delta_j}(E_j) \subset U_j.$ We notice that  $B_{\delta_j}(E_j) $ doesn't necessary contain $E_0$.

We prove the above by induction in $n$. For $n=1$, it is obvious. Assume the above statement holds for $n\leq p$. We consider the case $n=p+1$. We apply the induction for $a_1(E)\geq \cdots\geq a_p(E)$, i.e., there is an open set $U^p\subset B_\delta(E_0)$ and $0=i^p_0<i^p_1<\cdots<i^p_k=p$ such that
$$
a_{i_{m-1}+1}(E)=\cdots=a_{i_m}(E), \ \  1\leq m\leq k,
$$
for any $E\in U^p$.

Now we distinguish two cases:

{\bf Case I}:  $a_{p+1}(E)=a_p(E)$ for all $E\in U^p.$ Then choose $i^{p+1}_0=0$, $i^{p+1}_{1}=i^p_1$,..., $i^{p+1}_{k-1}=i^p_{k-1}$, $i^{p+1}_{k}=p+1$ and $U^{p+1}=U^p$.

{\bf Case II}:  $a_{p+1}(E')<a_p(E')$ for some $E'\in U^p.$ Then there is $U^{p+1}\subset U^{p}$ such that $a_{p+1}(E)<a_p(E)$ for all $E\in U^{p+1}$. We choose $i^{p+1}_0=0$, $i^{p+1}_{1}=i^p_1$,..., $i^{p+1}_{k}=i^p_{k}$ and $i^{p+1}_{k+1}=p+1$.

This completes the proof.
\end{pf}

By Lemma \ref{B2}, without loss of generality, we may assume there is a sequence of $E_n^j$ such that for each fixed $j$, there is $\delta_j>0$ such that
$$
\gamma_{n_{i-1}(E_n)}(E')>\gamma_{n_{i-1}(E_n)+1}(E')= \cdots=\gamma_{m(j)}(E')>\gamma_{m(j)+1}(E')=\cdots=\gamma_{n_{i}(E_n)}(E')> \gamma_{n_{i}(E_n)+1}(E')
$$
for $|E'-E_n^j|<\delta_j$.  By the definition of $n_i$,  there exists $s(E')\in\Z$ with
\begin{equation*}
n_{i-1}(E_n)= n_{s-1}(E'),    \quad  m(j)=n_{s}(E'),    \qquad  n_{i}(E_n)= n_{s+1}(E'),
\end{equation*}
such that
$$ \gamma_{n_{s-1}(E')}(E') > \gamma_{n_{s}(E')}(E') >\gamma_{n_{s+1}(E')}(E') .$$
Applying Lemma \ref{usel1}  to the turning points $\gamma_{m(j)}(E_n^j)$ and $\gamma_{n_i(E_n)}(E_n^j)$, we have
$$
\omega_{n_i(E^j_n)}^+(E^j_n)-\omega_{n_i(E^j_n)}^-(E^j_n)=  n_{s+1}(E_n^j)-n_s(E_n^j) .
$$
$$
\omega_{m(j)}^+(E^j_n)-\omega_{m(j)}^-(E^j_n)=n_s(E_n^j)-n_{s-1}(E_n^j).
$$

On the other hand, for any fixed $\e$ with
$-\gamma_{n_{i-1}(E_n)}(E_n)<\e<-\gamma_{n_i(E_n)}(E_n)$, the cocycle $(\alpha, A_{E_n}(\cdot+\e))$ is regular,  with acceleration
$$  \omega(\alpha, A_{E_n}(\cdot+ i \e))=  \omega_{n_i(E_n)}^-(E_n),$$ thus by Lemma \ref{regfor}, for $j$ sufficiently large,  $-\gamma_{n_{s-1}(E^j_{n})}(E^j_n) <\e<  -\gamma_{n_{s}(E^j_{n})}(E^j_n)$, such that  $(\alpha, A_{E^j_n}(\cdot+ i \e))$ is also regular, with acceleration $  \omega(\alpha, A_{E_n}(\cdot+ i \e))= \omega(\alpha, A_{E^j_n}(\cdot+ i \e))$, i.e.
 $$   \omega_{n_i(E_n)}^-(E_n) =\omega_{n_s(E^j_n)}^-(E^j_n)=\omega_{m(j)}^-(E^j_n). $$
 Similarly one can obtain
$$  \omega_{n_i(E_n)}^+(E_n) = \omega_{n_i(E^j_n)}^+(E^j_n).$$
 Consequently, by noting $\omega_{m(j)}^+(E^j_n) =\omega_{n_i(E^j_n)}^-(E^j_n)$, one has
\begin{align*}
\omega_{n_i(E_n)}^+(E_n)-\omega_{n_i(E_n)}^-(E_n)&=\omega_{n_i(E^j_n)}^+(E^j_n)-\omega_{n_i(E_n)}^-(E^j_n)+\omega_{m(j)}^+(E^j_n)-\omega_{m(j)}^-(E^j_n)\\
&= n_{s+1}(E_n^j)-n_s(E_n^j) +n_s(E_n^j) -n_{s-1}(E_n^j) = n_{i}(E_n)-n_{i-1}(E_n).
\end{align*}
\qed

\section{Proofs of the remaining results}\label{remain}
In this section, we give proofs of the remaining results in the introduction.\\
\noindent {\bf Proof of Theorem \ref{1general}:} Assume that
$$
\frac{\hat{L}_{k_1}(E)}{2\pi}<\frac{\hat{L}_{k_2}(E)}{2\pi}<\cdots< \frac{\hat{L}_{k_\ell}(E)}{2\pi}
$$
are the turning points of $L_\e(E)$ when $\e>0$, with the variation of the slope $\{k_i(E)-k_{i-1}(E)\}_{i=1}^\ell$ respectively where $k_0(E)=0$ and $k_\ell=m$. Thus we can express  $L_\e(E)$  as
\begin{align*}
L_{\e}(E)=\left\{
\begin{aligned} & L_0(E) &|\e|\in
\left[0,\frac{\hat{L}_{k_1}}{2\pi}\right]\\
&L_{\frac{\hat{L}_{k_i}}{2\pi}}(E)+2\pi k_{i}\left(|\e|-\frac{\hat{L}_{k_i}(E)}{2\pi}\right)
&|\e|\in \left(\frac{\hat{L}_{k_{i}}}{2\pi},\frac{\hat{L}_{k_{i+1}}}{2\pi}\right]\\
&L_{\frac{\hat{L}_{k_\ell}}{2\pi}}(E)+2\pi k_\ell\left(|\e|-\frac{\hat{L}_{k_\ell}(E)}{2\pi}\right) &|\e|\in \left(\frac{\hat{L}_{k_\ell}}{2\pi},h\right)
\end{aligned}\right.
\end{align*}
where $1\leq i\leq \ell-1$. Thus for any sufficiently small $\delta>0$,  and for any $\frac{\hat{L}_{k_\ell}}{2\pi}<h'<h$   we know that $(\alpha,A_E(\cdot+i\e))$ is regular if
$$
\e\in \begin{cases}
\left[\delta,\frac{\hat{L}_{k_1}}{2\pi}-\delta\right]\cup \left[\frac{\hat{L}_{k_\ell}}{2\pi}+\delta,h'\right]\cup \bigcup\limits_{i=1}^{\ell-1}\left[\frac{\hat{L}_{k_i}}{2\pi}+\delta,\frac{\hat{L}_{k_{i+1}}}{2\pi}-\delta\right]& \hat{L}_{k_1}>0\\
\left[\frac{\hat{L}_{k_\ell}}{2\pi}+\delta,h'\right]\cup \bigcup\limits_{i=1}^{\ell-1}\left[\frac{\hat{L}_{k_i}}{2\pi}+\delta,\frac{\hat{L}_{k_{i+1}}}{2\pi}-\delta\right]& \hat{L}_{k_1}=0
\end{cases}.
$$

Now we fix $E\in\R$, $h'<h$ and $V\in C^\omega_h(\T,\R).$ Let
$V^d(x)=\sum\limits_{k=-d}^d V_ke^{2\pi ikx},$ so we have
$$
\lim\limits_{d\rightarrow \infty}|V^d-V|_{h'}\rightarrow 0.
$$
By Lemma \ref{regfor}, there exists a neighborhood $I\times J$ of
$(V,\e)$, such that if $(V',\e')\in I\times J$, then $(\alpha,S_E^{V'}(\cdot+i\e'))$ is also regular, with acceleration
$$    \omega(\alpha,S_E^V(\cdot+i\e))= \omega(\alpha,S_E^{V'}(\cdot+i\e'))$$
where 
$$
S_E^V(x)=\begin{pmatrix}E-V(x)&-1\\ 1&0\end{pmatrix}.
$$
Consequently by a compactness argument,  for $d$ sufficiently large depending on $\delta,V$, and for any   $1\leq i\leq \ell-1$,    we have\\

\textbf{Case I: $\hat{L}_1>0$}
\begin{align*}
\omega(\alpha, S_E^{V^d}(\cdot+i\e))=\left\{
\begin{aligned} & 0 &|\e|\in
\left[\delta,\frac{\hat{L}_1}{2\pi}-\delta\right],\\
&k_{i}
&|\e|\in \left[\frac{\hat{L}_{i}}{2\pi}+\delta,\frac{\hat{L}_{i+1}}{2\pi}-\delta\right],\\
&k_m &|\e|\in \left[\frac{\hat{L}_m}{2\pi}+\delta,h'\right].
\end{aligned}\right.
\end{align*}

\textbf{Case II: $\hat{L}_1=0$}
\begin{align*}
\omega(\alpha, S_E^{V^d}(\cdot+i\e))=\left\{
\begin{aligned}
&k_{i}
&|\e|\in \left[\frac{\hat{L}_{i}}{2\pi}+\delta,\frac{\hat{L}_{i+1}}{2\pi}-\delta\right],\\
&k_m &|\e|\in \left[\frac{\hat{L}_m}{2\pi}+\delta,h'\right].
\end{aligned}\right.
\end{align*}

On the other hand, by Theorem \ref{1}, we have
\begin{align*}
L(\alpha, S_E^{V^d}(\cdot+i\e))=\left\{
\begin{aligned} & L(\alpha, S_E^{V^d}) &|\e|\in
\left[0,\frac{\hat{L}^d_{n_1}}{2\pi}\right]\\
&L\left(\alpha, S_E^{V^d}(\cdot+i\frac{\hat{L}^d_{n_s}}{2\pi})\right)+2\pi n_{s}\left(|\e|-\frac{\hat{L}^d_{n_s}(E)}{2\pi}\right)
&|\e|\in \left(\frac{\hat{L}^d_{n_{s}}}{2\pi},\frac{\hat{L}^d_{n_{s+1}}}{2\pi}\right]\\
&L\left(\alpha, S_E^{V^d}(\cdot+i\frac{\hat{L}^d_{n_\ell}}{2\pi})\right)+2\pi n_\ell\left(|\e|-\frac{\hat{L}^d_{n_\ell}(E)}{2\pi}\right) &|\e|\in \left(\frac{\hat{L}^d_{n_\ell}}{2\pi},h\right)
\end{aligned}\right.
\end{align*}
where $1\leq s\leq \ell-1$.  Hence for any $1\leq i \leq m$, there exists $s_1, s_2 \in \Z$ such that
\begin{eqnarray*}
n_{s_1}&=k_{i-1}, \ \ n_{s_2}&=k_{i},
\end{eqnarray*}
and furthermore, we have
$$\left|\frac{\hat{L}^d_{j}}{2\pi}-\frac{\hat{L}_{i}}{2\pi}\right|\leq 2\delta,\ \ n_{s_1}+1\leq j\leq n_{s_2}.$$
This actually  implies that
$$\left|\frac{\hat{L}^d_{j}}{2\pi}-\frac{\hat{L}_{i}}{2\pi}\right|\leq 2\delta,\ \ k_{i-1}+1\leq j\leq k_i.$$
Letting $\delta\rightarrow 0$, we get the result.\qed
\\

\noindent
{\bf Proof of Corollary \ref{3general}:} We distinguish two cases. If  $\hat{L}_1(E)>0$, then $L_\e(E)=L(E)$ for $|\e|\leq \hat{L}_1(E)$,  so  $\omega(E)=0$   by  definition. Otherwise, if  $\hat{L}_1(E)=0$, then by Theorem \ref{1general}, one has
$$L_\e(E)=L(E)+2\pi k_1(E)\e,    \qquad \e \in (0,\hat{L}_2(E)),$$
which implies  $
\omega(E)=k_1(E).
$\qed
\\

\noindent
{\bf Proof of Theorem \ref{2general}:}  We only prove (4), the other statements follow similarly. Note by \cite{avila}, $(\alpha, A_E)$ is subcritical, if
  and only if  $L(E)=0$ and $\omega(E)=0$. Then the result follows from
Corollary \ref{3general}. \qed \\

\noindent
{\bf Proof of Corollary \ref{4general}:}   It follows directly  from the definition of $h(E)$ and  Theorem \ref{1}.\qed
\\

\noindent
{\bf Proof of Corollary \ref{spec cha}:}   It follows directly  from Theorem \ref{jog} and (1) of Theorem \ref{2general}.\qed
\\

\noindent
{\bf Proof of Corollary \ref{localization}:}   It follows  from Corollary \ref{4general} and (1) of Theorem 1.2 in \cite{gyzh1}.\qed
\\


\noindent
{\bf Proof of Corollary \ref{6}:} It was proved in \cite{gyzh} (see Proposition 2.2 and Proposition 2.3) that if $(\alpha, \widehat{A}_{E})$ is almost reducible to some constant matrix $\widetilde{A}$, then
$$
\{\gamma_j\}_{j=1}^d=\{\ln|\lambda_j|\}_{j=1}^{d},
$$
where $\lambda_1, \cdots, \lambda_d$ are the eigenvalues of  $\widetilde{A}$ outside the unit circle, counting the multiplicity. Combining this fact with Theorem \ref{1}, we finish the proof. \qed
\\

\noindent
{\bf Proof of Corollary \ref{ghm1}:}  By Theorem \ref{thou}, we have
$$
L(E)=\gamma_1(E)+\gamma_2(E)+\ln|\lambda_2|.
$$
By Corollary \ref{3general}, $\omega(E)=2$ if and only if $\gamma_1(E)=\gamma_2(E)=0$. Thus the corollary follows.
\qed
\\

\noindent
{\bf Proof of Corollary \ref{ghm2}:}  Again by Theorem \ref{thou}, we have
$$
L(E)=\gamma_1(E)+\gamma_2(E)+\ln|\lambda_2|.
$$
Thus if $|\lambda_2|<1$, we must have $\gamma_1(E)>0$, then the results follow directly.
\qed

\section{Proof of Proposition \ref{mg}}\label{proofmg}
It suffices to prove
\be\label{gi}
\frac{\partial\omega_i}{\partial E}(E)=\tr\int_{\T}G_i(\theta)d\theta.
\ee
Once we have this, then the result follows from \eqref{reo} and Cauchy-Riemann equation.

To prove \eqref{gi}, we first need the following simple observation:
\begin{Lemma}\label{B1}
For $A\in C^1(\C,\rm{GL}_n(\C))$, we have
$$
\frac{\partial\ln \det{A(t)}}{\partial t}=\tr \frac{\partial A(t)}{\partial t} A^{-1}(t).
$$
\end{Lemma}
\begin{pf}
Let $A_{ij}$ be the $(i,j)$-th cofactor of the matrix $A.$ One can compute
\begin{align*}
\frac{\partial\ln \det{A(t)}}{\partial t}&=\frac{1}{\det{A(t)}}\frac{\partial \det{A(t)}}{\partial t}=\frac{1}{\det{A(t)}}\sum\limits_{i=1}^n\sum\limits_{j=1}^n\frac{\partial \det{A(t)}}{\partial a_{ij}} \frac{\partial a_{ij}(t)}{\partial t}\\
&=\frac{1}{\det{A(t)}}\sum\limits_{i=1}^n\sum\limits_{j=1}^nA_{ij}(t)\frac{\partial a_{ij}(t)}{\partial t}=\tr \frac{\partial A(t)}{\partial t} A^{-1}(t),
\end{align*}
thus the result follows.
\end{pf}

Consequently by Proposition \ref{f1} and Lemma \ref{B1}, we can compute
\begin{align*}
 \frac{\partial\omega_i}{\partial E}(E) =&\frac{\partial\int_{\T}\ln \det{M_-^i(\theta,E)}d\theta}{\partial E} \\ \nonumber
 =&\tr \int_\T \frac{\partial M^i_-(\theta,E)}{\partial E} (M^i_-(\theta,E))^{-1}d\theta\\ \nonumber
= &\tr\int_{\T} P_{[n_{i-1}+1,n_{i}]}\frac{\partial\widetilde M_-(T\theta,E)}{\partial E}\widetilde M^{-1}_-(T\theta,E)P_{[n_{i-1}+1,n_{i}]}d\theta. \nonumber
\end{align*}
Thus we need to compute $\frac{\partial\widetilde
  M_-(T\theta,E)}{\partial E}\widetilde M^{-1}_-(T\theta,E)$. Note
that by Proposition \ref{f1}, we have
\begin{equation}\label{tm}
\widetilde{M}_-(T\theta,E)=B^{-1}(\theta,E)M_-(T\theta,E)B(T\theta,E),
\end{equation}
therefore,  we have the following:
\begin{align*}
\begin{split}
&B^{-1}(\theta,E)\frac{\partial M_-(T\theta,E)}{\partial E} M^{-1}_-(T\theta,E)B(\theta,E)\\
=&B^{-1}(\theta,E)\frac{\partial \left (B(\theta,E)\widetilde M_-(T\theta,E) B^{-1}(T\theta,E)\right)}{\partial E} M^{-1}_-(T\theta,E)B(\theta,E)\\
=&B^{-1}(\theta,E)\frac{\partial B(\theta,E)}{\partial E}\widetilde M_-(T\theta,E) B^{-1}(T\theta,E) M^{-1}_-(T\theta,E)B(\theta,E)\\
&+B^{-1}(\theta,E)B(\theta,E)\frac{\partial \widetilde M_-(T\theta,E)}{\partial E} B^{-1}(T\theta,E) M^{-1}_-(T\theta,E)B(\theta,E)\\
&+B^{-1}(\theta,E)B(\theta,E)\widetilde M_-(T\theta,E)\frac{\partial B^{-1}(T\theta,E)}{\partial E}M^{-1}_-(T\theta,E)B(\theta,E)\\
=&\frac{\partial \widetilde M_-(T\theta,E)}{\partial E}\widetilde M^{-1}_-(T\theta,E) + E_2(\theta)
\end{split}
\end{align*}
where  we denote
\begin{equation}\label{e2}
E_2(\theta)=-\widetilde M_-(T\theta,E)B^{-1}(T\theta,E)\frac{\partial B(T\theta,E)}{\partial E}\widetilde M^{-1}_-(T\theta,E)+B^{-1}(\theta,E)\frac{\partial B(\theta,E)}{\partial E}.
\end{equation}

On the other hand, if we introduce the  auxiliary function
 $$F(\theta,E)=B^{-1}(\theta,E)G(\theta,E)\frac{\partial C^*M_-(\theta,E)}{ \partial E}B(\theta,E),$$
 we have the following observation:
\begin{Lemma}
 We have
\begin{align}\label{ft}
\begin{split}
& F(\theta,E)- \widetilde M_-(T\theta,E)F(T\theta,E)\widetilde M^{-1}_-(T\theta,E) \\
 =&\frac{\partial \widetilde M_-(T\theta,E)}{\partial E}\widetilde M^{-1}_-(T\theta,E)  -B^{-1}(\theta,E)G(\theta,E)B(\theta,E) + E_2(\theta).
\end{split}
\end{align}
\end{Lemma}
\begin{pf}
By \eqref{rce2}, we have
\begin{align*}
\frac{\partial C^*M_-(\theta,E)}{\partial E}=&-\frac{\partial CM_-^{-1}(T\theta,E)}{\partial E}-I_d=CM^{-1}_-(T\theta,E)\frac{\partial M_-(T\theta,E)}{\partial E}M^{-1}_-(T\theta,E)-I_d.
\end{align*}
Combining it with \eqref{g3},  one can compute
\begin{eqnarray*}
F(\theta,E)&=&B^{-1}(\theta,E)G(\theta,E)\left(CM^{-1}_-(T\theta,E)\frac{\partial M_-(T\theta,E)}{\partial E}M^{-1}_-(T\theta,E)-I_d \right)B(\theta,E)\\
&=&B^{-1}(\theta,E)G(\theta,E)CM^{-1}_-(T\theta,E)\frac{\partial M_-(T\theta,E)}{\partial E}M^{-1}_-(T\theta,E)B(\theta,E)\\
&&-B^{-1}(\theta,E)G(\theta,E)B(\theta,E) \\
&=& B^{-1}(\theta,E)M_-(T\theta,E)G(T\theta,E)  \frac{\partial C^* M_-(T\theta,E)}{\partial E}M^{-1}_-(T\theta,E)B(\theta,E)\\
 &&+ B^{-1}(\theta,E)\frac{\partial M_-(T\theta,E)}{\partial E}M^{-1}_-(T\theta,E)B(\theta,E)  -B^{-1}(\theta,E)G(\theta,E)B(\theta,E) \\
 &=& \widetilde M_-(T\theta,E)F(T\theta,E)\widetilde M^{-1}_-(T\theta,E) \\
  &&+ B^{-1}(\theta,E)\frac{\partial M_-(T\theta,E)}{\partial E}M^{-1}_-(T\theta,E)B(\theta,E)  -B^{-1}(\theta,E)G(\theta,E)B(\theta,E)
\end{eqnarray*}
where the last equality follows from \eqref{tm}. Thus the result follows.
\end{pf}

Note that $
\widetilde{M}_-(\theta,E)=\text{diag} \{M_-^1(\theta,E),M_-^2(\theta,E),\cdots, M_-^\ell(\theta,E)\}
$ is block diagonal. A direct computation shows that
\begin{align}\label{equ3}
&\tr\int_\T P_{[n_{i-1}+1,n_{i}]}E_2(\theta) P_{[n_{i-1}+1,n_{i}]}d\theta \\ \nonumber
=&\tr\int_T  M^i_-(T\theta,E)P_{[n_{i-1}+1,n_{i}]}B^{-1}(T\theta,E)\frac{\partial B(T\theta,E)}{\partial E}P_{[n_{i-1}+1,n_{i}]}(M^i_-(T\theta,E))^{-1}d\theta\\ \nonumber
&-\tr\int_T  P_{[n_{i-1}+1,n_{i}]}B^{-1}(\theta,E)\frac{ \partial B(\theta,E)}{\partial E}P_{[n_{i-1}+1,n_{i}]}d\theta\\ \nonumber
=&\tr\int_T  P_{[n_{i-1}+1,n_{i}]}B^{-1}(T\theta,E)\frac{\partial B(T\theta,E)}{\partial E}P_{[n_{i-1}+1,n_{i}]}d\theta\\ \nonumber
&-\tr\int_T  P_{[n_{i-1}+1,n_{i}]}B^{-1}(\theta,E)\frac{\partial B(\theta,E)}{\partial E}P_{[n_{i-1}+1,n_{i}]}d\theta=0.
\end{align}
The same argument shows that
\begin{align*} &\tr\int_\T P_{[n_{i-1}+1,n_{i}]}  \widetilde M_-(T\theta,E)F(T\theta,E)\widetilde M^{-1}_-(T\theta,E)P_{[n_{i-1}+1,n_{i}]}d\theta \\
&=  \tr\int_\T P_{[n_{i-1}+1,n_{i}]} F(\theta,E)P_{[n_{i-1}+1,n_{i}]}d\theta.
\end{align*}
Consequently by \eqref{ft}, we get the desired result.  \qed

\section{Proof of Lemma \ref{basic}}\label{proofbasic}
First we define the sequence
$$
u(n)=\begin{cases}
\frac{\sum\limits_{i=1}^m\phi_i(n)\Phi_{1,i}(q)}{a_d\det{\Phi(q)}} &\text{$n\geq q+1$}\\
-\frac{\sum\limits_{i=m+1}^{2d}\phi_i(n)\Phi_{1,i}(q)}{a_d\det{\Phi(q)}} &\text{$n\leq q$}
\end{cases}.
$$

For a fixed $n\in\Z$, one needs to consider the following cases:\\

\textbf{Case I:} If $n\geq q+d$ or $n\leq q-d$, then  it's straightforward to verify that $(L-EI)u(n)=0$.  \\

\textbf{Case II:} If $q+1\leq n< q+d$, we have
\begin{align*}
&(L-EI)u(n)=\sum\limits_{k=-d}^{q-n}a_ku(n+k)+\sum\limits_{k=q-n+1}^{d}a_ku(n+k)+(V(n)-E)u(n)\\
=&-\sum\limits_{k=-d}^{q-n}a_k\frac{\sum\limits_{i=m+1}^{2d}\phi_i(n+k)\Phi_{1,i}(q)}{a_d\det{\Phi(q)}} +\sum\limits_{k=q-n+1}^{d}a_k\frac{\sum\limits_{i=1}^m\phi_i(n+k)\Phi_{1,i}(q)}{a_d\det{\Phi(q)}}+(V(n)-E)\frac{\sum\limits_{i=1}^{m}\phi_i(n)\Phi_{1,i}(q)}{a_d\det{\Phi(q)}} \\
=&-\frac{\sum\limits_{k=-d}^{q-n}\sum\limits_{i=m+1}^{2d}a_k\phi_i(n+k)\Phi_{1,i}(q)}{a_d\det{\Phi(q)}} +\frac{\sum\limits_{i=1}^m(\sum\limits_{k=q-n+1}^{d}a_k\phi_i(n+k)+(V(n)-E)\phi_i(n))\Phi_{1,i}(q)}{a_d\det{\Phi(q)}}\\
=&-\frac{\sum\limits_{k=-d}^{q-n}\sum\limits_{i=m+1}^{2d}a_k\phi_i(n+k)\Phi_{1,i}(q)}{a_d\det{\Phi(q)}} -\frac{\sum\limits_{i=1}^m\sum\limits_{k=-d}^{q-n}a_k\phi_i(n+k)\Phi_{1,i}(q)}{a_d\det{\Phi(q)}} \\
=&-\frac{\sum\limits_{k=-d}^{q-n}a_k\sum\limits_{i=1}^{2d}\phi_i(n+k)\Phi_{1,i}(q)}{a_d\det{\Phi(q)}}
\end{align*}
where the penultimate equality follows from  the fact that $\phi_i$ is a solution of $Lu=Eu$.

On the other hand, by the assumption, we have $q-d+1\leq n+k\leq q.$
Thus  we always have $\sum\limits_{i=1}^{2d}\phi_i(n+k)\Phi_{1,i}(q)=0$. i. e., $(L-EI)u(n)=0$. \\

\textbf{Case III:}
If $q-d+1\leq n\leq q$, we have
\begin{align*}
&(L-EI)u(n)=\sum\limits_{k=-d}^{q-n}a_ku(n+k)+\sum\limits_{k=q-n+1}^{d}a_ku(n+k)+(V(n)-E)u(n)\\
=&-\sum\limits_{k=-d}^{q-n}a_k\frac{\sum\limits_{i=m+1}^{2d}\phi_i(n+k)\Phi_{1,i}(q)}{a_d\det{\Phi(q)}} +\sum\limits_{k=q-n+1}^{d}a_k\frac{\sum\limits_{i=1}^m\phi_i(n+k)\Phi_{1,i}(q)}{a_d\det{\Phi(q)}}-(V(n)-E)\frac{\sum\limits_{i=m+1}^{2d}\phi_i(n)\Phi_{1,i}(q)}{a_d\det{\Phi(q)}} \\
=&-\frac{\sum\limits_{i=m+1}^{2d}(\sum\limits_{k=-d}^{q-n}a_k\phi_i(n+k)+(V(n)-E)\phi_i(n))\Phi_{1,i}(q)}{a_d\det{\Phi(q)}} +\frac{\sum\limits_{i=1}^m\sum\limits_{k=q-n+1}^{d}a_k\phi_i(n+k)\Phi_{1,i}(q)}{a_d\det{\Phi(q)}} \\
=&\frac{\sum\limits_{i=m+1}^{2d}\sum\limits_{k=q-n+1}^{2d}a_k\phi_i(n+k)\Phi_{1,i}(q)}{a_d\det{\Phi(q)}} +\frac{\sum\limits_{i=1}^m\sum\limits_{k=q-n+1}^{d}a_k\phi_i(n+k)\Phi_{1,i}(q)}{a_d\det{\Phi(q)}} \\
=&\frac{\sum\limits_{k=q-n+1}^{d}a_k\sum\limits_{i=1}^{2d}\phi_i(n+k)\Phi_{1,i}(q)}{a_d\det{\Phi(q)}} .
\end{align*}
Note that if $n=q$, then $\frac{\sum\limits_{k=q-n+1}^{d}a_k\sum\limits_{i=1}^{2d}\phi_i(n+k)\Phi_{1,i}(q)}{a_d\det{\Phi(q)}}=1$, and by the assumption  if $n\leq q-1$, then $q+1\leq n+k\leq q+d-1$, thus $\frac{\sum\limits_{k=q-n+1}^{d}a_k\sum\limits_{i=1}^{2d}\phi_i(n+k)\Phi_{1,i}(q)}{a_d\det{\Phi(q)}}=0$.

Hence by the above discussions $(L-EI)u=\delta_q$, and it is obvious
that $u\in\ell^2(\Z)$, thus completing the proof.\qed

\appendix

\section{Proof of Theorem \ref{jog}}

We first prove the if part. If $(T,S_E^V)$ is uniformly hyperbolic, then by the definition, for any $x\in\Omega$, one can find two solutions $u_\pm(\cdot,x,E)\in \ell^2(\Z^\pm)$  obeying
$$
u_\pm(k-1,x,E)+u_\pm(k+1,x,E) + V(T^kx)u_\pm(k,x,E)=Eu_\pm(k,x,E).
$$
Therefore, by Proposition \ref{basic}, $(H_x-E)^{-1}$ exists and is
bounded for any $x\in \Omega$, so by Lemma \ref{notde}, $E\notin \Sigma$.

For the only if part we need the following result of Saker-Sell
\cite{sase}. The key is that the result works for complex valued potentials:

\begin{Lemma}\label{sasel}
If there are no non-trivial bounded solutions $u$ satisfying 
$H_x u=Eu$ for some $x$, then $(T,S_E^V)$ is  uniformly hyperbolic.
\end{Lemma}

Therefore, if  $(T,S_E^V)$ is not uniformly hyperbolic, then by Lemma \ref{sasel},  we can find a bounded vector $u$ such that $H_xu=Eu$ for some $x\in \Omega$. Consequently, we define $$
u_L(n)=
\begin{cases}
u(n)& |n|\leq L\\
0 & |n|>L
\end{cases}.
$$
A direct computation shows that
$$
\left\|(H_x-E)\frac{u_L}{\|u_L\|}\right\|^2=\frac{\|u_{L+1}\|^2-\|u_{L-1}\|^2}{\|u_{L}\|^2}\leq \frac{C}{\|u_{L}\|^2}.
$$
Note that we only need to consider the case $\|u_{L}\|^2\rightarrow\infty$, otherwise $E$ is an eigenvalue. In this case, $\frac{u_L}{\|u_L\|}$ is a Weyl sequence, hence $E\in \Sigma_x=\Sigma$.\qed

\section*{Acknowledgements}
This work was started   in 2015 when Q. Zhou was a  Visiting Assistant
Specialist at UCI, and completed  in 2020 when  L. Ge was a Visiting
Assistant Professor at UCI.  S. Jitomirskaya  was a 2020-21 Simons
fellow. Her work was also partially supported by NSF DMS-2052899,
DMS-2155211, and Simons 681675. She is also grateful to School of
Mathematics at Georgia Institute of Technology where she worked when
the manuscript was finalized.
J. You and Q. Zhou were partially supported by National Key R\&D Program of China (2020 YFA0713300) and Nankai Zhide Foundation.  J. You was also partially supported by NSFC grant (11871286). Q. Zhou was also supported by NSFC grant (12071232), the Science Fund for Distinguished Young Scholars of Tianjin (No. 19JCJQJC61300).

\end{document}